\documentclass[a4paper]{article}

\usepackage{lineno,hyperref}
\usepackage{amssymb}
\usepackage{amsmath}
\usepackage{amsfonts}
\usepackage{amsthm}
\modulolinenumbers[5]
\newtheorem{theorem}{Theorem}[section]
\newtheorem{lemma}[theorem]{Lemma}
\newtheorem{corollary}[theorem]{Corollary}
\newtheorem{proposition}[theorem]{Proposition}
\newtheorem{remark}{Remark}
\numberwithin{equation}{section}
\usepackage[utf8]{inputenc}
\usepackage[T1]{fontenc}

\begin{document}

\title{On asymptotic behaviour of solutions of the Dirac system
 and  applications to the Sturm--Liouville problem with a singular potential}
\date{}
\maketitle

\author{\textbf{{\L}ukasz Rzepnicki} \\ Faculty of Mathematics and Computer Science
Nicolaus Copernicus University
 Chopina 12/18, 87-100 Toruń, Poland\\ \verb"keleb@mat.umk.pl"\\
 
 \textbf{Alexander Gomilko} \\ Faculty of Mathematics and Computer Science
Nicolaus Copernicus University
 Chopina 12/18, 87-100 Toruń, Poland\\ \verb"alex@gomilko.com"
 }

\begin{abstract}
The main focus of this paper is
 the following matrix Cauchy problem for the Dirac system on the interval $[0,1]:$
\[
D'(x)+\left[\begin{array}{cc}
0 & \sigma_1(x)\\
\sigma_2(x) & 0
\end{array}
\right]
D(x)=i\mu\left[\begin{array}{cc}
1 & 0\\
0 &-1
\end{array}
\right]D(x),\quad D(0)=\left[\begin{array}{cc}
1 & 0\\
0 & 1
\end{array}\right],
\]
where %matrix $D=\{d_{ij}(x)\}_{i,j=1}^2$,
$\mu\in\mathbb{C}$ is a spectral parameter,
 and $\sigma_j\in L_2[0,1]$, $j=1,2$.
We propose a new approach for the study of asymptotic behaviour of its solutions as $\mu\to \infty$ and $|{\rm Im}\,\mu|\le d$.
As an application,  we obtain new, sharp
asymptotic formulas for eigenfunctions
of  Sturm--Liouville operators with singular potentials.
\end{abstract}
\noindent
\textbf{keyword:}\\
Sturm--Liouville problem, singular potential,
Dirac system, Sturm--Liouville operator\\
\verb"MSC"[2010] 34L20,  34E10

\section{Introduction}

The motivation to study the Dirac systems came from our interest in the
spectral problem for Sturm--Liouville equation with a singular potential,
in particular the asymptotic behaviour with respect to a spectral parameter of its solutions.
We found a way to transform the spectral problem into a special form of the perturbed Dirac system.
The known results related to asymptotic behaviour of the solutions of the Dirac system
were not precise enough for our aims. Therefore we attracted our attention to the Dirac system itself.
In this way we obtained new results about Dirac systems and gained a new tool
to spectral analysis of Sturm--Liouville problems.

The paper thus is naturally divided into two parts: one is devoted to
Dirac systems and the other one is focused on the applications of the results to spectral problems for Sturm--Liouville equations with singular potentials.
We start with the presentation of the background concerning Dirac systems.

Consider  the following matrix Cauchy problem for the Dirac system:
\begin{equation}\label{CaychyInt}
D'(x)+J(x)D(x)=
i\mu J_0 D(x),\quad x\in [0,1],\quad D(0)=I,
\end{equation}
where
%matrices
\begin{equation}\label{matrixA}
J_0=\left[\begin{array}{cc}
1 & 0\\
0 &-1
\end{array}
\right],\quad
J(x)=\left[\begin{array}{cc}
0 & \sigma_1(x)\\
\sigma_2(x) & 0
\end{array}
\right],\quad
I:=\left[\begin{array}{cc}
1 & 0 \\
0 & 1
\end{array}\right],
\end{equation}
$\mu\in \mathbb{C}$ is a spectral parameter, and  $\sigma_j\in L_2[0,1]$, $j=1,2$, are complex-valued functions.
By a solution of \eqref{CaychyInt}, as usual, we understand a
matrix $D$ with entries from
the space of absolutely continuous on $[0,1]$ functions
(i.e. from the Sobolev space $W_1^1[0,1]$)
satisfying  \eqref{CaychyInt} for a.e. $x\in [0,1]$. In fact here both conditions yield that $D$ has entries from $W_2^1[0,1]$

In recent years the study of asymptotic behavior  of  fundamental solutions to
Dirac  systems has attracted a lot of attention, and we may mention e.g.
\cite{Alb1}, \cite{Amirov},  \cite{Siegl},
\cite{DM0},
\cite{DM},
\cite{DM2},
\cite{Levi},
\cite{LMal},
\cite{Sad1}, \cite{S}, and \cite{SSDirac}, as samples.
In particular, there is a variety of  asymptotic formulas for the solutions to \eqref{CaychyInt}
in different settings and under a wide range of assumptions on the potentials.
However, most of those formulas  contain just
a leading term.
In this paper, by a new and comparatively simple argument, we find
sharp asymptotic for solutions to Dirac systems, thus generalizing
most of similar relations from the literature.

More specifically, we study the asymptotic behaviour of solutions $D(x)=D(x,\mu)$ to the Cauchy
problem \eqref{CaychyInt}
as $\mu\to \infty$, $\mu\in  P_d$, where
\[
P_d:=\{\mu\in {\rm C}:\,|{\rm Im}\,\mu|\le d\}.
\]
The main results related to this issue are placed in Section 2 in Lemma \ref{asD} and Corollary \ref{Cor3.2}.

Our results can be analyzed wider since there exists the connection between the Dirac system presented in this paper and more general formulations of
the problem, often studied in literature.
Recall that any solution $V=\{v_1,v_2\}^T$, $v_j\in W_1^1[0,1],$  of the next Cauchy problem for the Dirac system:
\begin{equation}\label{System1}
V'(x)+J(x)V(x)=i\mu J_0V(x),\quad x\in [0,1],
\end{equation}
\[
\quad V(0)=\{c_1,c_2\}^T\in \mathbb{C}^2,
\]
is of the form
\[
V(x)=D(x)V(0),\quad x\in [0,1].
\]

Note that a  system formally more general than \eqref{System1}:
\begin{equation}\label{System2}
Y'(x)+P(x)Y(x)=i\mu J_0Y(x),\quad P(x)=\left[
\begin{array}{cc}
p_{11}(x) & p_{12}(x)\\
p_{12}(x) & p_{22}(x)\end{array}
\right],
\end{equation}
where $x\in [0,1]$ and $Y=\{y_1,y_2\}^T$,
$p_{jj}\in L_1[0,1]$, $p_{ij}\in L_2[0,1]$, $i,j=1,2$, $i\not=j$,
can be reduced to \eqref{System1} by a simple transformation given in
\cite{Sad1}. Indeed,
setting
\[
r_j(x)=\int_0^x p_{jj}(t)\,dx,\quad j=1,2,
\]
it is easy to see that  $Y$ is a solution of \eqref{System2} if and only if a vector-function
\[
V(x)=\left[\begin{array}{cc}
e^{r_1(x)} & 0\\
0 & e^{r_2(x)}\end{array}\right]Y(x)
\]
satisfies  \eqref{System1} with
\[
\sigma_1(x)=p_{12}(x)e^{r_1(x)-r_2(x)},\quad
\sigma_2(x)=p_{21}(x)e^{r_2(x)-r_1(x)}.
\]
On the other hand, the classical Dirac system
(see. e. g.  \cite[Ch. VII, $\S\,1$]{Levi})
\begin{equation}\label{System3}
BZ'(x)+Q(x)Z(x)
=\mu Z(x),\quad x\in [0,1],
\end{equation}
where
\[
B=\left[\begin{array}{cc}
0 & 1\\
-1 &  0
\end{array}
\right],\quad  Q(x)=\left[
\begin{array}{cc}
q_{11}(x) & q_{12}(x)\\
q_{12}(x) & q_{22}(x)\end{array}
\right],\quad q_{ij}\in L_2[0,1],
\]
can be transformed to \eqref{System2} by writing
\[
Z(x)=UY(x),\quad
U=\left[\begin{array}{cc}
1 & -i\\
-i &  1
\end{array}\right].
\]

The idea of our method, originating from \cite[Ch. 1, $\S\,2$, Problem 5]{March},
relies on  a special
integral representation
for  $D(x)=D(x,\mu)$.
The representation has the form
\begin{align}\label{D0}
D(x)=e^{xA_\mu}+\int_0^xe^{(x-2t)A_\mu}
[J(t)+Q(x,t)]\,dt,
\end{align}
where
\begin{align}\label{amu}
A_\mu:=i\mu J_0=\left[\begin{array}{cc}
i\mu & 0\\
0 &-i\mu
\end{array}
\right],
\end{align}
and $Q$ is a matrix function
continuous on
\begin{equation}\label{DDelta}
\Delta:=\{(x,t)\in \mathbb{R}^2: \;  0\leq t \leq x\leq1\},
\end{equation}
and, moreover, $Q$ is
a unique solution in the space of continuous functions on $\Delta$ of the
integral equation
\begin{align}\label{Qtilde0}
Q(x,t)=
\int_0^{x-t}J(t+\xi)J(\xi)d\xi
-\int_0^{x-t}J(t+\xi)Q(t+\xi,\xi)d\xi.
\end{align}
It is crucial to observe that
the kernel $Q$
does not depend on $\mu$. More details can be found in Lemma \ref{jednor}.
This approach allows us to obtain much sharper asymptotic formulas for solutions to
Dirac systems than those in the literature. As a consequence,
we get fine asymptotics for solutions  to perturbed Dirac systems of the form
\begin{equation}\label{tildeS}
\tilde{D}'(x)+J(x)\tilde{D}(x)=i\mu J_0 \tilde{D}(x)+\frac{P(x)}{\mu}\tilde{D}(x),
\quad x\in [0,1],\quad \tilde{D}(0)=I,
\end{equation}
where $\mu\not=0$ and the entries of $P$ are in $L_1[0,1]$.
Note that such an approach can be extended to a more general case
when $\sigma_j\in L_p[0,1]$, $j=1,2$, $p\ge 1$. Perturbed Dirac systems are analyzed
in Section 3, where the main result is Theorem \ref{ArCa}.

%??????We do not consider this case
%in present article in view of
%application our result to differential equation with singular potential,
%where as we show sufficient consider
%However, we restrict ourselves the case $\sigma_j\in L_2[0,1]$, $j=1,2$.

It should be mentioned that the representation \eqref{D0} is not quite new.
In \cite{Alb1} the authors
considered a matrix Dirac equation
\begin{equation}\label{ee}
BU'(x)+\left[\begin{array}{cc}
q_1(x) & q_2(x) \\
q_2(x)  & -q_1(x)\end{array}\right]U(x)=\lambda U(x),\qquad x\in [0,1],\quad U(0)=I,
\end{equation}
where $q_j$, $j=1,2$ are real-valued functions in $L_p[0,1]$, $p\in [1,\infty)$.
By \cite[Theorem 2.1]{Alb1}
one has
\begin{equation}\label{Lit1}
U(x,\lambda)=e^{-\lambda x B}+
\int_0^x e^{-\lambda(x-2s)} P(x,s)\,ds,
\end{equation}
where for every $x\in [0,1]$ the entries $p_{ij}(x,\cdot)$, $i,j=1,2$ of $P(x,\cdot)$, belong to $L_p[0,1]$ and the mapping
$x\mapsto p_{ij}(x,\cdot)$ is continuous on $[0,1]$.
The representation \eqref{Lit1} was derived in \cite{Alb1}
via the method of successive approximations
by solving the integral equation
\[
U(x,\lambda)=e^{-\lambda x B}+
\int_0^x e^{-\lambda(x-t)B} B Q(t)U(t,\lambda)\,dt,\quad x\in [0,1],\quad Q=\{q_{ij}\}_{i,j=1}^2,
\]
and using the property $e^{-\lambda x B}Q(t)=Q(t)e^{\lambda x B}$, $x,t\in [0,1]$.
In the Hilbert space setting, when  $p=2,$ we follow a different strategy for obtaining   \eqref{D0}.
Namely, we deduce an  integral equation for $Q$
and investigate this equation directly. As a result,
we obtain more detailed information on $Q$ and this allows to derive more precise asymptotic
formulas for the solutions of the Dirac system (see Lemma \ref{asD}).

Note that in \cite{S} and \cite{SSDirac}
the Dirac system
\eqref{System3} was studied in a situation when
$q_{ij}$, $i,j=1,2$ are complex-valued functions in $L_p[0,1]$, $p\in [1,\infty)$.
By means of  Pr\"ufer substitution \cite{Hartman}, in a rather complicated manner,
  ``short'' asymptotic formulas (as $\mu\to\infty$, $\mu\in P_d$)
were obtained in \cite{S} and \cite{SSDirac}
for  solutions $Z$ of %the system
\eqref{System3} with  the initial conditions
$Z(0)=\{1,0\}^T$ and $Z(0)=\{0,1\}^T$.
These formulas (with $p=2$) correspond
to our Corollary \ref{Cor3.2}.

We show in Section 4 and 5 that our results apply to the study of
asymptotic behaviour of the system of fundamental solutions to the differential equation
\begin{equation}\label{eqS}
l(y)+\lambda y=0,\qquad l(y):=y''+q(x)y,\quad x\in [0,1],
\end{equation}
where a potential $q$ belongs to $W^{-1}_2[0,1]$. In other words,
\[
q(x)=\sigma'(x),\qquad \sigma\in L_2[0,1],
\]
where the derivative is understood in the sense of distributions.
Such a class of singular potentials includes
Dirac $\delta$-type and Coulomb
$1/x$-type interactions, which are frequently
used in quantum
mechanics and mathematical
physics. For a comprehensive treatment of  physical models with potentials from
negative Sobolev spaces  we refer to \cite{Kappeler}.
A short history and different approaches to the study of singular
Sturm--Liouville operators can be found in \cite{SS}.

Following the regularization method described in
\cite{Atk} and \cite{SS}, given $y$
we formally introduce a quasi-derivative $y^{[1]}$ of $y$:
\[
y^{[1]}(x):=y'(x)+\sigma(x)y(x),
\]
so that \eqref{eqS} can be rewritten in the form
\begin{equation}\label{Seq}
(y^{[1]})'(x)-\sigma(x)y^{[1]}(x)+\sigma^2(x)y(x)+\lambda y(x)=0,\qquad x\in [0,1].
\end{equation}

Our main idea is to transform the equation
\eqref{Seq} into an appropriate
perturbed Dirac system
\eqref{tildeS}.
We show that if  $\mu\not=0$ then the transformation
\begin{equation}\label{trans}
\left(\begin{array}{c}
y\\
y^{[1]}\end{array}\right)=\left[\begin{array}{cc}
1 & 1\\
i\mu & -i\mu\end{array}\right]
\left(\begin{array}{c}
v_1\\
v_2\end{array}\right)
\end{equation}
``recasts''  \eqref{Seq}
as the following perturbed Dirac system for $V=\{v_1,v_2\}^T$, $x\in [0,1]:$
\[
V'(x)+\left[\begin{array}{cc}
0 & \sigma(x)\\
\sigma(x) & 0\end{array}\right]V(x)=i\mu J_0 V(x)+
\frac{i\sigma^2(x)}{2}\left[\begin{array}{cc}
1 & 1\\
-1 & -1\end{array}\right]V(x).
\]
In this way, using our results
on asymptotic of Cauchy's matrix for \eqref{tildeS},
 we  obtain  very precise information on
asymptotic behavior of fundamental solutions to \eqref{eqS} and
on asymptotic distribution of the eigenvalues and eigenfunctions
to the spectral problem
\begin{equation}\label{specy0}
l(y)+\lambda y=0,\qquad
y(0)=0, \qquad  y(1)=0.
\end{equation}

The main results Theorem \ref{YSa} and Corollary \ref{fund} concerning the asymptotic behavior of fundamental solutions to \eqref{eqS}
are placed in Section 4.

In Section 5 we concentrate on spectral problem \eqref{specy0}.
The problem \eqref{specy0} was widely studied in the literature. Basic facts concerning the asymptotics of its eigenvalues
and eigenfunctions can be found in \cite{Hryn2} and \cite{Hryn}. The authors of \cite{Hryn2} and \cite{Hryn} employed
the  method of transformation operators
and used a factorized form of $l$ to obtain asymptotic formulas with remainder terms expressed in terms of  Fourier sine coefficients.
Their results were precise enough to study the inverse spectral problem for \eqref{specy0}, see, for instance, \cite{Alb1} and \cite{Alb2}.
More detailed formulas, stated in Theorem \ref{evalue},  were obtained in \cite{S0}, \cite{SS} and \cite{Sav2006}.
The approach of the latter papers relied on the Pr\"{u}fer substitution,  and it allowed one to derive
asymptotic formulas for eigenfunctions $y_n$ of \eqref{specy0}
in the form
\begin{align}\label{asst}
\pi n y_n(x)=y_{0,n}(x)+\psi_{n}(x), \ \ n\geq1,
\end{align} where $y_{0,n}$
is a known function (the explicit formulation of it is given by \eqref{tezafw200}),
and  remainders $\psi_{1,n}$ satisfying
\begin{equation}\label{ff}
\sup_{x\in[0,1]}\sum_{n=1}^\infty |\psi_{n}(x)|<\infty.
\end{equation}
This method was further developed in a number of subsequent articles, see e.g.
% by various authors in the series of papers on the
%spectral problem \eqref{specy0} with different types of initial conditions or the topic of equiconvergence.
%we refer the reader to
\cite{Vladyk}, \cite{Shveik2}, \cite{Shveik3}, and \cite{Shveik1}.
Our main  result
%, Theorem \ref{wlasne},
generalizes and sharpens  Theorem \ref{evalue} essentially.
In particular, we provide similar asymptotic formulas to \eqref{asst} i.e.
\begin{align}\label{asstt}
\pi n y_n(x)=y_{0,n}(x)+y_{1,n}(x)+\tilde{\psi}_{n}(x), \ \ n\geq 1,
\end{align}
where $y_{1,n}$ are known functions (see \eqref{tezafw}) and
remainder terms  $\tilde{\psi}_{n}$
satisfy a strengthened form  of \eqref{ff}:
\[
\sum_{n=1}^\infty\sup_{x\in[0,1]}| \tilde{\psi}_{n}(x)|<\infty.
\]
%a property much weaker than \eqref{ff}.
%holds.

Note that our method of reduction to the perturbed Dirac system applies in a setting
more general than that of this paper.
For example (see Remark \ref{remarc}), a transformation similar to
\eqref{trans}
allows one to treat  the equations of the form
\[
(a(x)y'(x))'+q(x)y(x)+\lambda c(x) y(x)=0,\qquad x\in [0,1],
\]
where $q=u',$ $u\in L_2[0,1],$
and
%the coefficients $a$ and $c$ satisfy the next assumptions:
\[
a\in W_2^1[0,1],\quad c\in W_2^1[0,1],
\quad a(x)>0,\quad c(x)>0,\quad x\in [0,1].
\]

The paper is organized as follows. Section 2 is devoted to the Dirac system \eqref{CaychyInt}. We
derive there a special representation and establish asymptotic formulas of the solutions.
In Section 3 we use the facts from the previous section to analyze the asymptotic behaviour of the
solutions of the perturbed Dirac system.
Section 4 is focused on the application of the results related to the Dirac system to investigation of Sturm--Liouville
equations with singular potentials, in particular we present there how to transform a Sturm--Liouville equation into the Dirac system of a special form.
The Dirichlet spectral problem for Sturm--Liouville equation with singular potential is studied in Section 5. For the clarity of exposition some part of auxiliary and technical results necessary in the proofs of Sections 2 and 3
are placed in the Appendix.

In the paper we adopt the following convention: when the domain of a function is not given, it is assumed
to be $[0,1]$. For example, we use the notation
$L_p:=L_p[0,1]$.

\section{Matrix Cauchy problems for Dirac system}

This section is devoted to the matrix Cauchy problem
\begin{equation}\label{v120}
D'(x)+J(x) D(x)=A_\mu D(x),\quad
D(0)=I,
%:=\left[\begin{array}{cc}
%1 & 0\\
%0 & 1
%\end{array}\right],
\ \ x\in [0,1],
\end{equation}
where  $A_\mu=i\mu J_0$ (see \eqref{amu})
and $\sigma_j\in L_2[0,1]$, $j=1,2$.

We obtain an important
integral representation
for the solutions of \eqref{v120}, given  by \eqref{D0}
and \eqref{Qtilde0}.

First, we introduce a necessary notation. Denote
the space of continuous functions
on $\Delta$ with the supremum norm $\|\cdot\|_{C(\Delta)}$ by $C(\Delta)$.
Define
\begin{equation}\label{coetants}
l_0:=\max\{\|\sigma_1\|_{L_1},\|\sigma_2\|_{L_1}\},\quad
l:=\|\sigma_1\|_{L_1}\cdot\|\sigma_2\|_{L_1},\quad
l_1:=\|\sigma_1\|_{L_1}+\|\sigma_2\|_{L_1},
\end{equation}
and
\begin{equation}\label{constants1}
\tilde{l}_0:=\max\{\|\sigma_1\|_{L_2},\|\sigma_2\|_{L_2}\},\quad
\tilde{l}:=\|\sigma_1\|_{L_2}\cdot\|\sigma_2\|_{L_2},
\quad
l_2=\|\sigma_1\|^2_{L_2}+\|\sigma_2\|^2_{L_2}.
\end{equation}
Moreover, let
\begin{equation}\label{NfunC}
\sigma_0(x):=|\sigma_1(x)|+|\sigma_2(x)|\in L_2[0,1].
\end{equation}

Let $X$ be a Banach space, then  $M(X)$ will stand for  the Banach
space of $2\times 2$ matrices with entries
from $X$ and the norm
\[
\|Q\|_{M(X)}:=\sum_{k,j=1}^2 \|Q_{jk}\|_X,\quad
Q=[Q_{jk}]_{j,k=1}^2.
\]

In our analysis of $D$ we will need several auxiliary results addressing  the two integrals in \eqref{Qtilde0}.
If $J_0$ and $J$ are defined by
\eqref{matrixA},
then
\begin{equation}\label{Note1}
J_0^2=I,\quad J_0 J(x)+J(x)J_0=0,\qquad \mbox{a.e.}\;x\in [0,1].
\end{equation}
Moreover, let  the matrix function $\tilde{J}$ be given by
\begin{equation}\label{ksigmazeroo}
\tilde{J}(x,t):=
\int_0^{x-t}J(t+\xi)J(\xi)d\xi=
\int_t^x J(s)J(s-t)\,ds,\quad(x,t)\in \Delta.
\end{equation}
Note that
\[
\tilde{J}(x,t)=\left(\begin{array}{cc}
\tilde{\sigma}_1(x,t) & 0\\
0 & \tilde{\sigma}_2(x,t)\end{array}\right),
\]
where
\begin{equation}\label{defSS}
\tilde{\sigma}_1(x,t):=\int_0^{x-t}\sigma_1(t+\xi)\sigma_2(\xi)d\xi,\quad
%\int_t^x \sigma_1(s)\sigma_2(s-t)\,ds,\quad(x,t)\in \Delta,
%\]
%\[
\tilde{\sigma}_2(x,t):=\int_0^{x-t}\sigma_2(t+\xi)\sigma_1(\xi)d\xi.
%\int_t^x \sigma_2(s)\sigma_1(s-t)\,ds,
\end{equation}

Using  Cauchy--Schwarz inequality, for  $(x,t)\in \Delta$ and
$(x+\epsilon,t+\delta)\in \Delta$, we have
\begin{align*}
|\tilde{\sigma}_1(x+\epsilon,t+\delta)&-\tilde{\sigma}_1(x,t)|\\&=
\Big|\int_0^{x-t+\epsilon-\delta}\sigma_1(t+\delta+s)\sigma_2(s)
ds-\int_0^{x-t}\sigma_1(t+s)\sigma_2(s)ds\Big|
\\
&\leq \left(\int_{x-t}^{x-t+\epsilon-\delta}|\sigma_2(s)|^2ds
\right)^{1/2}\cdot \left(\int_{x+\delta}^{x+\epsilon}|\sigma_1(s)|^2ds
\right)^{1/2}
\\
&+\|\sigma_2\|_{L_2}
\left(\int_t^x|\sigma_1(\delta+s)-
\sigma_1(s)|^2ds\right)^{1/2}.
\end{align*}
Clearly, a similar estimate holds  for $\tilde{\sigma}_2$ as well.
Thus $\tilde{\sigma}_j\in C(\Delta), j=1,2,$ and
\begin{equation}\label{L22CY}
\|\tilde{\sigma}_j\|_{C(\Delta)}\le \tilde{l}, \qquad j=1,2.
\end{equation}

To treat the other integral from \eqref{Qtilde0},
define  a linear
operator $T_{\sigma}$ on $C(\Delta)$  by
\begin{align}\label{T}
(T_{\sigma}f)(x,t)=\int_0^{x-t}\sigma(t+\xi)f(t+\xi,\xi)d\xi=
\int_t^x\sigma(s)f(s,s-t)ds, %\quad f\in C(\Delta),
\end{align}
where $\sigma\in L_1[0,1]$.
Let $(x,t)\in \Delta$
and $\epsilon\in \mathbb{R}$, $\delta\in \mathbb{R}$
so that $(x+\epsilon,t+\delta)\in \Delta$.
Then
\begin{align*}
|(T_\sigma f)(x+\epsilon,t+\delta)&-(T_\sigma f)(x,t)|
\\
&\le \int_{t+\delta}^{x+\epsilon}\sigma(s)|f(s,s-t-\delta)
-f(s,s-t)|ds
\\
&+\left|\int_{t+\delta}^{x+\epsilon}\sigma(s)f(s,s-t)\,ds-
\int_t^x \sigma(s)f(s,s-t)\,ds\right|
\\
&\le \|\sigma\|_{L_1}
\sup_{s\in [0,1]}
\sup_{\tau_1,\tau_2\in[0,s],\,|\tau_1-\tau_2|\le \delta}\,
|f(s,\tau_1)-f(s,\tau_2)|
\\
&+2\|f\|_{C(\Delta)}
\sup_{0<b-a<\max\{\epsilon,\delta\}}\int_a^b |\sigma(s)|\,ds.
\end{align*}
Therefore, due to the fact that $f$ is uniformly continuous  on $\Delta$,
$T_\sigma$ is bounded on $C(\Delta)$.
Moreover,
\begin{equation}\label{MoreA}
|(T_\sigma f)(x,t)|\le \|f\|_{C(\Delta)}\int_0^x|\sigma(t)|\,dt\le \|\sigma\|_{L_1}\|f\|_{C(\Delta)},\quad (x,t)\in \Delta,\quad f\in C(\Delta).
\end{equation}

Next, define bounded linear operators
$T_{12}$ and $T_{21}$ on $C(\Delta)$ by
\[
T_{12}:=T_{\sigma_1} T_{\sigma_2},\quad
T_{21}:=T_{\sigma_2} T_{\sigma_1},
\]
and remark that
\begin{equation}\label{Sthat}
(T_{kj}f)(x,t)=\int_t^x \sigma_k(s)
\left(\int_{s-t}^s \sigma_j(\tau)f(\tau,\tau-s+t)\,d\tau\right)\,ds,\quad k,j=1,2,\quad k\not=j.
\end{equation}
The operators $T_{kj}$, $k,j=1,2$, $k\not=j$, enjoy the following estimate.
\begin{lemma}\label{lemT}
One has
\begin{align}\label{Tn}
\|T_{kj}^nf\|_{C(\Delta)}
\le
\frac{l^n}
{n!}\|f\|_{C(\Delta)},\qquad f\in C(\Delta),\quad n\in \mathbb{N},\quad k,j=1,2,\quad k\not=j.
\end{align}
\end{lemma}

\begin{proof} Consider the operator $T_{12}$.
Define $\eta\in C[0,1]$ by
\[
\eta(x):=\int_0^x |\sigma_1(s)|\left(\int_0^s |\sigma_2(\tau)|\,d\tau\right)\,ds,\qquad x\in [0,1].
\]
It suffices to prove that
for all $(x,t)\in \Delta$ and $n=1,2,\dots,$
\begin{equation}\label{100}
|(T_{12}^nf)(x,t)|
\le \frac{\|f\|_{C(\Delta)}}{n!}
\eta^n(x),\qquad
 f\in C(\Delta).
\end{equation}
If $n=1$, then the estimate \eqref{100}
follows directly from  \eqref{Sthat}.
Arguing by induction, suppose that \eqref{100} holds for some $n\in \mathbb{N}$.
Then for $(x,t)\in \Delta$ and $f\in C(\Delta)$ we have
\begin{align*}
|(T_{12}^{n+1}f)(x,t)|&\le \int_t^x
|\sigma_1(s)|\int_{s-t}^s |\sigma_2(\tau)||(T_{12}^nf)(\tau,\tau-s+t)|\,d\tau\,ds
\\
&\le \frac{\|f\|_{C(\Delta)}}{n!}
\int_0^x |\sigma_1(s)| \int_{s-t}^s |\sigma_2(\tau)| \eta^n(\tau)\,d\tau\,ds
\\
&\le \frac{\|f\|_{C(\Delta)}}{n!}
\int_0^x |\sigma_1(s)| \int_0^s |\sigma_2(\tau)|\,d\tau\,\eta^n (s)\,ds
\\
&=\frac{\|f\|_{C(\Delta)}}{n!}
\int_0^x \eta^n(s)\,d\eta(s)
=\frac{\|f\|_{C(\Delta)}}{(n+1)!}\eta^{n+1}(x).
\end{align*}
Therefore \eqref{100} and then \eqref{Tn} hold true.
\end{proof}

Now we are ready to establish  a first crucial property of the solutions to \eqref{v120}.
The proof of the next statement
follows an idea from \cite[Ch. 1, $\S 24$]{March}.

\begin{lemma}\label{jednor}
If $\sigma_j\in L_2[0,1]$, $j=1,2$ and $\mu \in P_d$, then the unique solution $D=D(x,\mu)$ of \eqref{v120} %-\eqref{Dzero}
can be represented as
\begin{align}\label{D}
D(x,\mu)=e^{xA_\mu}+\int_0^xe^{(x-2t)A_\mu}[J(t)+Q(x,t)]dt,
\end{align}
where
$Q\in M(C(\Delta))$ is the unique solution of the integral equation
\begin{align}\label{Qtilde}
Q(x,t)=
\tilde{J}(x,t)+
\int_0^{x-t}J(t+\xi)Q(t+\xi,\xi)d\xi,
\end{align}
where $\tilde{J}\in M(C(\Delta))$  is given by \eqref{ksigmazeroo}.
Moreover,
\begin{equation}\label{L1Est}
\|Q\|_{M(C(\Delta))}\le
2\tilde{l}(1+l_0)e^l,
\end{equation}
and
\begin{equation}\label{LLl}
\|D\|_{M(C[0,1])}
\leq e^{d}a, \quad a:=1+l_1+2\tilde{l}(1+l_0)e^l.
\end{equation}
\end{lemma}

\begin{proof}
The uniqueness of solutions comes from \cite[Thm. 1.2.1]{Zettl}. We look for  solutions
of \eqref{v120} in the form
\begin{align}\label{wzor1}
D(x,\mu)=e^{xA_\mu}U(x,\mu),
\quad U(0,\mu)=I.
\end{align}
The relation
\begin{equation}\label{relAA}
J(x)e^{xA_\mu}=e^{-xA_\mu}J(x),\quad \mbox{a. e.}\;\;x\in [0,1]
\end{equation}
implies that $U$ satisfies the Cauchy problem
\[
U'(x,\mu)+e^{-2x A_\mu}J(x) U(x,\mu)=0,\qquad x\in [0,1],\quad U(0,\mu)=I,
\]
which is equivalent to the integral equation
\begin{align}\label{wzor15}
U(x,\mu)=I-\int_0^xe^{-2tA_\mu}J(t)U(t,\mu)dt,\qquad x\in [0,1].
\end{align}
Now, we look for solutions of
\eqref{wzor15}
in the form
\begin{align}\label{wzor2}
U(x,\mu)=I+\int_0^xe^{-2tA_\mu}Q_0(x,t)dt,\quad
\end{align}
 where $Q_0\in M(L_2(\Delta))$ does not depend on $\mu$.
 Substituting \eqref{wzor2} into \eqref{wzor15},
we obtain
\begin{align*}
\int_0^x e^{-2tA_\mu}Q_0(x,t)\,dt&=
-\int_0^x e^{-2tA_\mu}J(t)\,dt\\
&-\int_0^x e^{-2tA_\mu}J(t)
\int_0^t e^{-2s A_\mu} Q_0(t,s)\,ds\,dt.
\end{align*}
Then, by \eqref{Note1}, we have
\begin{align*}
\int_0^xe^{-2tA_\mu}J(t)
\int_0^t
e^{-2s A_\mu} &Q_0(t,s)ds\,dt
=\int_0^xe^{-2tA_\mu}
\int_0^t
e^{2s A_\mu} J(t)Q_0(t,s)ds\,dt
\\
&=\int_0^xe^{-2tA_\mu}\int_0^{x-t}
J(t+\xi)Q_0(t+\xi,\xi)d\xi dt,
\end{align*}
hence
\[
\int_0^xe^{-2tA_\mu}Q_0(x,t)\,dt=
-\int_0^xe^{-2tA_\mu}\left(J(t)
+\int_0^{x-t}
J(t+\xi) Q_0(t+\xi,\xi)d\xi\right)\,dt
\]
for all $x\in [0,1].$ Therefore, $U$ is a solution of \eqref{wzor15}
if and only if $Q_0\in M(L_2(\Delta))$ is a solution of
\begin{align}\label{prev}
Q_0(x,t)=-J(t)-\int_0^{x-t}J(t+\xi)Q_0(t+\xi,\xi)d\xi.
\end{align}
Next,  setting
\[
Q_0(x,t)=-J(t)+Q(x,t),\qquad (x,t)\in \Delta,
\]
and using \eqref{Qtilde}, we infer that  $Q$ satisfies \eqref{prev}.
The latter equation can be written in an operator form
\begin{align}\label{Ttilde}
Q=\tilde{J}+\tilde{T} Q,\quad
\tilde{T}=
-\left[\begin{array}{cc}
0 & T_{\sigma_1}\\
T_{\sigma_2} & 0
\end{array}\right],
\end{align}
where the operators $T_{\sigma_1}$ and $T_{\sigma_2},$  defined by \eqref{T}, are  linear and bounded  on $C(\Delta)$.
In particular, by \eqref{MoreA}, we have
\[
\|\tilde{T}F\|_{M(C(\Delta))}
\le
\max\{\|\sigma_1\|_{L_1},\|\sigma_2\|_{L_1}\}
\|F\|_{M(C(\Delta))},\quad F\in M(C(\Delta)).
\]
A more detailed analysis of $\tilde{T}$
can be found in the Appendix.

Next observe that
\[
\tilde{T}^{2n}=
\left[\begin{array}{cc}
T_{12}^n & 0\\
0 & T_{21}^n
\end{array}\right],\quad n\in {\rm N},
\]
therefore, by \eqref{Tn},
\[
\|\tilde{T}^{2n}F\|_{M(C(\Delta))}
\le
\frac{l^n}{n!}\|F\|_{M(C(\Delta))},\quad F\in M(C(\Delta)).
\]
So
\eqref{Qtilde} has a unique solution
$Q\in C(\Delta)$ of the form
\begin{equation}\label{Pred}
Q=\sum_{n=0}^\infty \tilde{T}^n
\tilde{J}=\sum_{n=0}^\infty \tilde{T}^{2n}(I+\tilde{T})\tilde{J},
\end{equation}
and moreover
\begin{equation}\label{Pred11}
\quad \|Q\|_{M(C(\Delta))}
\le (1+l_0)e^l\|\tilde{J}\|_{M(C(\Delta))}.
\end{equation}
Thus, \eqref{Pred} and \eqref{L22CY} imply \eqref{L1Est}.

Finally, in view of \eqref{D} and \eqref{L1Est}, we obtain
\begin{equation*}
\|D\|_{M(C[0,1])}\leq e^{d} \Big(1+l_1+\|Q(x,t)\|_{M(C(\Delta))}\big),\quad \mu\in P_d,
\end{equation*}
and  \eqref{LLl} follows.
\end{proof}

By means of simple calculations, one obtains
the following  representation of the matrix $Q$ defined by \eqref{Pred}.

\begin{corollary}\label{Skalar}
If $Q$ is given by
\eqref{Pred}, then
\begin{align}
Q=\begin{bmatrix}
Z_1 & -T_{\sigma_1} Z_2\\
-T_{\sigma_2} Z_1 & Z_2
\end{bmatrix},\quad
Z_1=\tilde{\sigma}_1+ \sum_{k=1}^\infty T_{12}^{k}\tilde{\sigma}_1,\quad
Z_2=\tilde{\sigma}_2+ \sum_{k=1}^\infty T_{21}^{k}\tilde{\sigma}_2.
\end{align}
Moreover,
\begin{equation}\label{HCC}
\|Z_j\|_{C(\Delta)}\leq \tilde{l}\exp{l},
\end{equation}
\[
\|T_{\sigma_j} Z_l\|_{C(\Delta)}\leq \|\sigma_j\|_{L_1}\|Z_l\|_{C(\Delta)}, j=1,2.
\]
\end{corollary}

%Note that
%\[
%T_{\sigma_2} Z_1=Z_2T_{\sigma_2},\quad T_{\sigma_1} Z_2=Z_1T_{\sigma_1}.
%\]

Going back to the formulas \eqref{D} and
\eqref{Pred} for
$D=D(x,\mu)$, $x\in [0,1]$, $\mu\in P_d$,
note that \begin{align}
D(x,\mu)&=e^{xA_\mu}+\int_0^xe^{(x-2t)A_\mu}[J(t)+Q(x,t)]dt
\nonumber\\
&=e^{xA_\mu}+\int_0^xe^{(x-2t)A_\mu}J(t)dt
+\int_0^xe^{(x-2t)A_\mu}\tilde{J}(x,t)dt
\nonumber\\
&+\int_0^xe^{(x-2t)A_\mu}(\tilde{T}\tilde{J})(x,t)dt\nonumber
\\&+\int_0^xe^{(x-2t)A_\mu}\sum_{n=2}^\infty (\tilde{T}^{n}\tilde{J})(x,t)dt.\label{DDD}
\end{align}

Now, we are ready to investigate
the asymptotic behaviour of  $D(x,\mu)$ defined by \eqref{DDD}
for $\mu \in P_d$ and fixed $d>0$.

In what follows we will use different types of
estimates for reminders.
For fixed $\sigma_j\in L_2$, $j=1,2$, and  $\mu\in \mathbb{C}$ define
\begin{align}\label{ga0}
\gamma(\mu)&:=\sum_{j=1}^2\left(\Big\|\int_0^xe^{-2i\mu t}\sigma_j(t)dt\Big\|_{L_2}+\Big\|\int_0^xe^{2i\mu t}\sigma_j(t)dt\Big\|_{L_2}
\right),
\end{align}
\begin{align}\label{ga10}
\gamma_0(x,\mu)&:=\sum_{j=1}^2\left(\Big|\int_0^xe^{-2i\mu t}\sigma_j(t)dt\Big|+\Big|\int_0^xe^{2i\mu t}\sigma_j(t)dt\Big|\right),\quad x\in [0,1].
\end{align}
and
\begin{equation}\label{gamma12}
\gamma_1(\mu):=\int_0^1 \sigma_0(s)\gamma_0^2(s,\mu)\,ds,\qquad
\gamma_2(\mu):=l_2\gamma^2(\mu)+l_1\gamma_1(\mu).
\end{equation}
It is easy to see that if $\mu\in P_d$, then
\begin{equation}\label{gS1}
\|\gamma_0(x,\mu)\|_{L_2}\le \gamma(\mu),\quad
\gamma(\mu)\le 2e^{2d}l_1,\quad
\gamma_0(x,\mu)\le 2e^{2d}l_1,\quad x\in [0,1],
\end{equation}
and
\begin{equation}\label{gamma2-1}
\gamma_2(\mu)\le
4e^{4d}l_1^2(l_2+l_1^2),\quad \gamma_2(\mu)\le 2l_1e^{2d}(l_2+2l_1e^{2d}\|\sigma_0\|_{L_2})\gamma(\mu),
\end{equation}

For the sequel, define
\begin{equation}\label{NotNN}
N(x,t):=(\tilde{J}+\tilde{T}\tilde{J})(x,t)\in C(\Delta),
\end{equation}
and observe  that
\[
N(x,t)= \left[\begin{array}{cc}
\tilde{\sigma}_1(x,t) & -(T_{\sigma_1}\tilde{\sigma}_2)(x,t)\\
-(T_{\sigma_2}\tilde{\sigma}_1)(x,t) &  \tilde{\sigma}_2(x,t)
\end{array}
\right]\quad (x,t)\in \Delta.
\]

\begin{lemma}\label{asD}
Let  $\sigma_j\in L_2$, $j=1,2$.
If $D(x,\mu)$ is a solution of \eqref{v120} then
\begin{align}\label{t1}
D(x,\mu)&=e^{x A_\mu} + D^{(0)}(x,\mu)+D^{(1)}(x,\mu),
\end{align}
where
\[
D^{(0)}(x,\mu)=\int_0^x e^{(x-2t)A_\mu} J(t)\,dt
+\int_0^x
e^{(x-2t)A_\mu} N(x,t)\,dt,
\]
and for all $\mu\in P_d$ and $x \in [0,1],$
\[
\|D^{(1)}(x,\mu)\|_{M(C[0,1]}\le 2e^{5d}\exp{(e^{2d}l_1)}\gamma_2(\mu).
\]
\end{lemma}

\begin{proof}
By \eqref{DDD} we have
\[
D^{(1)}(x,\mu)=\int_0^xe^{(x-2t)A_\mu}\sum_{n=2}^\infty (\tilde{T}^n\tilde{J})(x,t)dt.
\]
and using  the inequality \eqref{EstIm3} proved in the Appendix, we infer that
\begin{align*}
\|D^{(1)}(x,\mu)\|_{M(C[0,1])}
&\le \sum_{n=2}^\infty \left\|\int_0^xe^{(x-2t)A_\mu}(\tilde{T}^n\tilde{J})(x,t)dt\right\|_{M(C[0,1])}
\\
&\le 2 e^d \gamma_2(\mu)\sum_{n=2}^\infty \frac{e^{2nd}l_1^{n-2}}{(n-2)!}=
2e^{5d}\exp{(e^{2d}l_1)}\gamma_2(\mu),
\end{align*}
for all $x\in [0,1]$ and $\mu\in P_d$.
\end{proof}

The above lemma leads to asymptotic formulas for $D(x,\mu)$, $\mu\in P_d$.
Let us start with  several
simple observations.
%First of all, remark that
%\begin{align*}
%\|D^{(0)}(x,\mu)\|_{M(C([0,1]))}&\le e^d(\|J\|_{M(L_1[0,1])}
%+\|\tilde{J}\|_{M(C(\Delta))}
%+\|\tilde{T}\tilde{J}\|_{M(C(\Delta))})
%\\
%&\le e^d(l_1+(1+l_1)\|\tilde{J}\|_{M(C(\Delta))})\\
%&\le
%e^d(l_1+2(1+l_1)\tilde{l}).
%\end{align*}
Note that from
\[
\left\|\int_0^x e^{(x-2t)A_\mu} J(t)\,dt\right\|_{M(\mathbb{C})}\le
e^d\gamma_0(x,\mu),\quad x\in [0,1],
\]
\eqref{EstIm0}, and \eqref{EstIm100} it follows that
\[
\left\|D^{(0)}(x,\mu)\right\|_{M({\rm C})}\le
e^d\gamma_0(x,\mu)+ 2e^{3d} (1+\tilde{l}_0e^{2d})\tilde{l}_0 \gamma(\mu),\quad x\in [0,1].
\]
Moreover, by \eqref{EstIm0}, \eqref{EstIm100}, we have
\[
\left\|\int_0^x e^{(x-2t)A_\mu} N(x,t)\,dt\right\|_{M(C[0,1]}\le
2e^{3d}\tilde{l}_09!+e^{2d}\tilde{l}_0)\gamma(\mu).
\]
Combining  these inequalities with
 Lemma \ref{asD} and the estimates from
 \eqref{gamma2-1}, we obtain the following representations for $D$.

\begin{corollary}\label{Cor3.2}
For every  $d>0$ there exist  $a_j=a_j(d,\|\sigma_1\|_{L_2},\|\sigma_2\|_{L_2})$, $j=0,1,2$ such that
for all $x\in [0,1]$ and $\mu\in P_d,$
\begin{equation}\label{t222}
D(x,\mu)=e^{x A_\mu}+R(x,\mu),
\end{equation}
where
\[
\|R(x,\mu)\|_{M(C[0,1])}\le a_0,\quad
\|R(x,\mu)\|_{M(\rm C)}\le
 a_1(\gamma(\mu)+\gamma_0(x,\mu)),\quad x\in [0,1].
\]
Moreover,
\begin{equation}\label{LCdop}
D(x,\mu)=e^{tA_\mu}
+D_{0}(x,\mu)+R_0(x,\mu),
\end{equation}
where
\[
D_{0}(x,\mu):=\int_0^x e^{(x-2t)A_\mu}J(t)\,dt
\]
and
\[
\|R_0(x,\mu)\|_{M(C[0,1]}\le
 a_2(\gamma(\mu)+\gamma_1(\mu)),\quad x\in [0,1].
\]
\end{corollary}

\begin{remark}
Note that the explicit formula for $D_0$ is the following
\[
D_{0}(x,\mu)=\left[
\begin{array}{cc}
0 & q_1(x,\mu) \\
q_2(x,\mu) & 0
\end{array}
\right],
\]
\[
q_1(x,\mu):=\int_0^x e^{i\mu (x-2t)} \sigma_1(t)\,dt,\quad
q_2(x,\mu):=\int_0^x e^{-i\mu(x-2t)} \sigma_2(t)\,dt.\]
\end{remark}

\section{Perturbed Dirac system}

Let us consider the matrix Cauchy problem
\begin{equation}\label{CPr1}
\tilde{D}'(x)+J(x)\tilde{D}(x)=A_\mu \tilde{D}(x)+
\frac{P(x)}{\mu}\tilde{D}(x),\quad x\in [0,1],\quad
\tilde{D}(0)=I,
\end{equation}
where the entries of the matrix $P(x)=[p_{k,j}(x)]_{k,j=1}^2$ are complex-valued functions from $L_1[0,1].$
Let  $\mu\not=0$ and let $D$ be a solution of \eqref{v120}.
Then the solutions of
\eqref{CPr1} also satisfy the integral   equation
\begin{equation}\label{IntAm}
\tilde{D}=D+\frac{1}{\mu}\mathcal{A}\tilde{D},
\end{equation}
where the linear operator $\mathcal{A}=\mathcal{A}(\mu)$ acting on $M(C[0,1])$
is given by
\begin{align}\label{operatorA}
(\mathcal{A}\tilde{D})(x)=D(x)(\mathcal{A}_0\tilde{D})(x),\quad
(\mathcal{A}_0\tilde{D})(x):=\int_0^x D^{-1}(t) P(t) \tilde{D}(t)\,dt,\quad x\in [0,1].
\end{align}

Note that by Liouville's theorem,
$$
\det D(x)=\det D(0)=1, \qquad x\in [0,1],
$$
hence
\begin{equation}\label{Dminus}
D^{-1}(x)=-\tilde{J}_0[D(x)]^T\tilde{J}_0,\quad
\tilde{J}_0:=\left[\begin{array}{cc}
0 & 1\\
-1 & 0
\end{array}
\right], \qquad  x\in [0,1].
\end{equation}
As a consequence, for every $x\in [0,1],$
\begin{equation}\label{DDE}
\|D(x)\|_{M(\mathbb{C})}=\|D^{-1}(x)\|_{M(\mathbb{C})}.
%\quad
%\|D\|_{M(C[0,1])}=\|D^{-1}\|_{M(C[0,1])}.
\end{equation}
The identity \eqref{DDE} and the estimate \eqref{LLl}
imply that
%the norm of the operator $\mathcal{A}(\mu)$ in
%${M(C[0,1])}$ admits the estimate
\begin{equation}\label{OperA}
\|\mathcal{A}(\mu)\|_{M(C[0,1])}\le \|D\|^2_{M(C[0,1])}
\|P(t)\|_{M(L_1[0,1])}\le be^{2d},\quad \mu\in P_d,
\end{equation}
where
\[
b:=a^2\|P\|_{M(L_1[0,1])},
\]
and a constant $a$ is defined in  \eqref{LLl}.

The following statement relates the solutions of the perturbed matrix Cauchy problem
\eqref{CPr1}
to the solutions of the  Cauchy problem
\eqref{v120}.

\begin{corollary}\label{mainCor}
For every $\mu\in P_d$,
$|\mu|>b e^{2d}$, the solution $\tilde{D}(x,\mu)$
of  \eqref{CPr1}
admits the representation
\begin{align}\label{szereg}
\tilde{D}(x,\mu)=D(x,\mu)+\sum_{k=1}^{\infty}
\mu^{-k}\mathcal{A}^k(\mu)D(x,\mu),
\end{align}
where the series converges in $M(C[0,1])$.
Moreover, for $\mu\in P_d$, and $|\mu|>2be^{2d},$
\begin{align}\label{omu}
\|\tilde{D}-D\|_{M(C[0,1])}\leq \frac{a_2}{|\mu|},\quad a_2=2 a b e^{3d}
\end{align}
and
\begin{equation}\label{duzeV2}
\tilde{D}(x,\mu)=D(x,\mu)
+\frac{1}{\mu}\mathcal{A}(\mu)D(x,\mu)+
\tilde{D}_1(x,\mu),
\end{equation}
where
\begin{align}\label{reminder}
\|\tilde{D}_1(x,\mu)\|_{M(C[0,1])}\le
\sum_{k=2}^{\infty}
|\mu|^{-k}\|\mathcal{A}^k(\mu)
D(x,\mu)\|_{M(C[0,1])}
\le \frac{a_3}{|\mu|^2},
\end{align}
and $a_3=2 a b^2e^{5d}$.
\end{corollary}

Now we focus on the representation \eqref{duzeV2}. In particular we
investigate the behaviour of the second term $1/\mu\mathcal{A}(\mu)D(x,\mu)$
in \eqref{duzeV2}.
Due to \eqref{Dminus}
and \eqref{t222}, we have
\begin{align}
D^{-1}(t)P(t)D(t)&=-\tilde{J}_0D^T(t)\tilde{J}_0 P(t)D(t)\nonumber\\
&=
-\tilde{J}_0[e^{tA_\mu}+R^T(t,\mu)]\tilde{J}_0 P(t)
[e^{tA_\mu}+R(t,\mu)]
\nonumber\\
&=e^{-tA_\mu} P(t)
e^{tA\mu}+R_1(t,\mu),\label{Dminuss}
\end{align}
where for almost all $t\in [0,1],$
\begin{align*}
\|R_1(t,\mu)\|_{M(\mathbb{C})}&\le
(\|R\|^2_{M(\mathbb{C})}+
2e^{d}
\|R\|_{M(\mathbb{C})})\|P(t)\|_{M(\mathbb{C})}
\\
&\le a_4(\gamma(\mu)+\gamma_0(t,\mu))\|P(t)\|_{M(\mathbb{C})},
\end{align*}
with a constant $a_4=a_4(d,\|\sigma_1\|_{L_2},\|\sigma_2\|_{L_2})$.
Note that
\[
e^{-tA_\mu} P(t)
e^{tA\mu}
=\left[\begin{array}{cc}
p_{11}(t) & e^{-2 i\mu t}p_{12}(t) \\
e^{2i\mu t}p_{21}(t) & p_{22}(t)
\end{array}\right].
\]
Furthermore, we have
\begin{align*}
\left\|
\int_0^x R_1(t,\mu)\,dt\right\|_{M( C[0,1])}
&\le  a_4(\|P\|_{M(L_1[0,1])}
\gamma(\mu)+
k_P(\mu)),
\end{align*}
where
\begin{align}\label{k}
k_P(\mu):=\int_0^1 \|P(t)\|_{M(\mathbb{C})}\gamma_0(t,\mu)\,dt.
\end{align}

Using \eqref{LCdop}, we obtain that
\begin{align*}
\mathcal{A}(\mu)D(x,\mu)&=
\left[e^{A_\mu x}
+D_0(x,\mu)+R_0(x,\mu)\right]
\int_0^x
e^{-tA_\mu} P(t)
e^{tA\mu}
\,dt
\\
&+D(x,\mu)\int_0^x R_1(t,\mu)\,dt
\\
&=e^{xA_\mu}\int_0^x
e^{-tA_\mu} P(t)
e^{tA\mu}
\,dt
\\&+D_0(x,\mu)
\int_0^x
e^{-tA_\mu} P(t)
e^{tA\mu}
\,dt
+R_2(x,\mu),
\end{align*}
where
\begin{align*}
R_2(x,\mu)&=R_0(x,\mu)
\int_0^x e^{-tA_\mu} P(t)
e^{tA\mu}\,dt+D(x,\mu)\int_0^x R_1(t,\mu)\,dt
\end{align*}
and
\[
\|R_2(x,\mu)\|_{M( C[0,1])}\le
a_5\Big(\|P\|_{M(L_1[0,1])}\big(\gamma(\mu)+\gamma_1(\mu)\big)+k_P(\mu)\Big),\quad \mu\in P_d,
\]
with a constant $a_5=a_5(d,\|\sigma_1\|_{L_2}, \|\sigma_2\|_{L_2}).$

In this way, combining the observations following Corollary \ref{mainCor} with \eqref{duzeV2}, we derive
the following assertion yielding sharp asymptotics for  $\tilde{D}$.

\begin{theorem}\label{ArCa}
If $\mu\in P_d$, $|\mu|>2be^{2d}$, then
\[
\tilde{D}(x,\mu)=\mathcal{R}(x,\mu)
+\mathcal{R}_0(x,\mu),\quad x\in [0,1],
\]
where $\mathcal{R}$ is given by
\begin{align*}
\mathcal{R}(x,\mu)&:=e^{xA_\mu}+D^{(0)}(x,\mu)
+\frac{1}{\mu}e^{xA_\mu}\int_0^x
e^{-tA_\mu} P(t)
e^{tA\mu}\,dt
\\&+\frac{1}{\mu}
D_0(x,\mu)
\int_0^x e^{-tA_\mu} P(t)
e^{tA\mu}\,dt,
\end{align*}
where
\[
\mathcal{R}_0(x,\mu)=
D^{(1)}(x,\mu)+\frac{R_2(x,\mu)}{\mu}
+\tilde{D}_1(x,\mu),
\]
with $D^{(1)}$ given by \eqref{t1}
and the reminder satisfies
\begin{align*}
\|\mathcal{R}_0(x,\mu)
\|_{M(\mathbb{C})}&
\le a_6(\gamma^2(\mu)+\gamma_1(\mu)+|\mu|^{-1}k_P(\mu)+|\mu|^{-2}),
\end{align*}
with a constant $a_6=a_6(d,\|\sigma_1\|_{L_2},\|\sigma_2\|_{L_2},\|P\|_{M(L_1)})$.
\end{theorem}

\section{Sturm--Liouville equations with singular potentials}

In this section we consider a Sturm--Liouville equation
\begin{equation}\label{SE}
y''(x)+q(x)y(x)+\mu^2  y(x)=0,\quad x\in [0,1],\quad
\end{equation}
where the potential $q$ is a complex-valued distribution
of the first order,  i.e.
\[
q=\sigma',\quad \sigma\in L_2[0,1],
\]
and $\mu\in \mathbb{C}$ is a spectral parameter.
Following the regularization method, we
introduce the quasi-derivative of $y\in W_1^1[0,1]$ as
\[
y^{[1]}(x):=y'(x)+\sigma(x)y(x),\qquad
y'(x)=y^{[1]}(x)-\sigma(x)y(x),
\]
and rewrite  \eqref{SE} as
\begin{equation}\label{SE10}
(y^{[1]}(x))'-\sigma(x)y^{[1]}(x)+\sigma^2(x)y(x)
+\mu^2 y(x)=0,\qquad x\in [0,1].
\end{equation}
We say that  $y$ is a solution of \eqref{SE10} if
\[
y\in \mathcal{D}:=\{y\in W_1^1[0,1],\quad
y^{[1]}(x)\in W_1^1[0,1]\},
\]
and  \eqref{SE10} is satisfied  for a.e. $x\in [0,1]$.

Note that \eqref{SE10} can be written in a matrix form
\begin{equation}\label{Znew}
L\left(\begin{array} {c}
y \\
y^{[1]}\end{array}
\right):=
\frac{d}{dx}
\left(\begin{array} {c}
y \\
y^{[1]}\end{array}
\right)+
\left[\begin{array}{cc}
\sigma & -1\\
\sigma^2+\mu^2 & -\sigma\end{array}\right]\left(\begin{array} {c}
y \\
y^{[1]}\end{array}
\right)=0.
\end{equation}
Let us now relate the  equation  \eqref{SE10}
to the perturbed system \eqref{CPr1}.
For $\mu\not=0$ define
\[
S_1(\mu):=\left[\begin{array}{cc}
1 &  0\\
0 & i\mu
\end{array}\right]S_0,
\quad
S_0:=\left[\begin{array}{cc}
1 & 1 \\
1 & -1
\end{array}\right],\quad
S_1^{-1}(\mu)=
\frac{1}{2i\mu}\left[
\begin{array}{cc}
i\mu & -1\\
i\mu & -1
\end{array}\right].
\]
If
\begin{equation}\label{Change1}
\left(\begin{array} {c}
y \\
y^{[1]}\end{array}
\right)=S_1(\mu)
\left(\begin{array} {c}
v_1 \\
v_2\end{array}
\right),\quad
\left(\begin{array} {c}
v_1 \\
v_2\end{array}
\right)=S_1^{-1}(\mu)
\left(\begin{array} {c}
y \\
y^{[1]}\end{array}
\right),
\end{equation}
then by \eqref{Znew},
\[
L\left(\begin{array} {c}
y \\
y^{[1]}\end{array}
\right)=S_1(\mu)L_0\left(\begin{array} {c}
v_1 \\
v_2\end{array}
\right),
\]
where
\[
L_0\left(\begin{array} {c}
v_1 \\
v_2\end{array}
\right):=\left(\begin{array} {c}
v_1 \\
v_2\end{array}
\right)'+
\frac{1}{2i\mu}
\left[\begin{array}{cc}
2\mu^2+\sigma^2 & 2i\mu \sigma+\sigma^2\\
2i\mu \sigma-\sigma^2 & -2\mu^2-\sigma^2\end{array}\right]\left(\begin{array} {c}
v_1 \\
v_2\end{array}
\right).
\]
Summarizing, we showed that $y$ is a solution of
\eqref{SE10}
if and only if $V=\{v_1,v_2\}^T$, defined by
\eqref{Change1}, is a solution of the system
\begin{equation}\label{CaychyInt0}
V'(x)+\sigma(x)\left[\begin{matrix}
0 & 1\\
1 & 0\end{matrix}\right]V(x)
=i\mu J_0 V(x)+\frac{i\sigma^2(x)}{2\mu}\left[\begin{matrix}
1 & 1\\
-1 & -1\end{matrix}\right]V(x),\quad x\in [0,1],%\quad V(0)=\vec{c}\in \mathbb{C}^2,
\end{equation}
and this is precisely \eqref{CPr1}
with
\[
J(x)=\left[\begin{array}{cc}
0 &  \sigma(x)\\
\sigma(x) & 0\end{array}\right],\quad \sigma\in L_2[0,1],
\]
\begin{equation}\label{rfun0}
P(x)=\frac{i}{2}\left[\begin{array}{cc}
\tau(x) &  \tau(x)\\
-\tau(x) & -\tau(x)\end{array}\right],\quad \tau=\sigma^2\in L_1[0,1].
\end{equation}
Note that in this case we have $\|\tau\|_{L_1}=\|\sigma\|_{L_2}^2$.

As a result, if $y$ is a solution of \eqref{SE10} with the initial conditions
\begin{equation}\label{initS}
y(0)=c_1,\quad y^{[1]}(0)=c_2,
\end{equation}
and $\mu\not=0$ then we can
apply Theorem \ref{ArCa}
to \eqref{CaychyInt0}. In this case
\[
 e^{-tA_\mu} P(t)
e^{tA\mu}\,dt=\frac{i}{2}\left[\begin{array}{cc}
1 & e^{-2i\mu t}\\
-e^{-2i\mu t} & 1
\end{array}\right],
\]
and
\[
D_0(x,\mu)=\int_0^x \left[\begin{array}{cc}
0 & e^{i\mu (x-2t)}\\
e^{-i\mu (x-2t)} & 0\end{array}\right]\,\sigma(t)\,dt,
\]
\[
D^{(0)}(x,\mu)=D_0(x,\mu)+\int_0^x e^{(x-2t)A_\mu}
\left[
\begin{array}{cc}
\tilde{\sigma}(x,t) & -(T_\sigma\tilde{\sigma})(x,t)\\
 -(T_\sigma\tilde{\sigma})(x,t) & \tilde{\sigma}(x,t)\end{array}\right]\,dt,
\]
where
\[
\tilde{\sigma}(x,t)=\int_0^{x-t}\sigma(t+\xi)\sigma(\xi)\,d\xi,
\]
\begin{align}\label{nsigma}
N_\sigma(x,t):=(T_\sigma\tilde{\sigma})(x,t)=\int_0^{x-t} \sigma(t+\xi)\tilde{\sigma}(t+\xi,\xi)\,d\xi.
\end{align}
Theorem \ref{ArCa} and the transformation \eqref{Change1} yield the identity
\[
\left(\begin{array}{c}
y \\
y^{[1]}
\end{array}\right)
=S_1(\mu)\left[\mathcal{R}(x,\mu)+
\mathcal{R}_0(x,\mu)\right]S_1^{-1}(\mu)\left(\begin{array}{c}
c_1\\ c_2\end{array}\right).
\]
Thus we proved the following statement in which and
in what follows we will write
\begin{equation}\label{rhoA}
\rho(\mu):=\gamma^2(\mu)+ \gamma_1(\mu)+|\mu|^{-1}
k_{\sigma^2}(\mu)]
+|\mu|^{-2},
\end{equation}
where all necessary definitions are given by \eqref{ga0}, \eqref{ga10} \eqref{gamma12} and $\eqref{k}$ with
$\sigma_1=\sigma_2=\sigma$.

\begin{theorem}\label{YSa}
Let $d>0$ and
\begin{equation}\label{domain}
\Omega_d(\sigma):=\{\mu\in P_d: \, |\mu|>2be^{2d}\},\;\;
b=\|\sigma\|_{L_2}^2\left(
1+\|\sigma\|_{L_1}+
e^{\|\sigma\|_{L_1}}\|\sigma\|_{L_2}^2\right)^2.
\end{equation}
If $y=y(x,\mu)$, $\mu\in \Omega_d(\sigma)$, is a solution of
\eqref{SE10}
 with the initial conditions
\eqref{initS}, then
\[
\left(\begin{array}{c}
y\\
(i\mu)^{-1}y^{[1]}\\
\end{array}
\right)=W(x,\mu)\vec{c}_0(\mu)+
W_0(x,\mu)\vec{c}_0(\mu),
\]
where
\[
\vec{c}_0(\mu)=S_1^{-1}(\mu)\left(\begin{array}{c}
c_1\\ c_2\end{array}\right)=
\frac{1}{2}\left(\begin{array}{c}
c_1+\frac{c_2}{i\mu}\\
c_1-\frac{c_2}{i\mu}\end{array}\right),
\]
and
\begin{align*}
W(x,\mu)
&=\left[\begin{array}{cc}
e^{i\mu x} & e^{-i\mu x}\\
e^{i\mu x} & -e^{-i\mu x}\end{array}
\right]\\&+
\left[\begin{array}{cc}
d_1(x,\mu)+d_2(x,-\mu) &
d_2(x,\mu)+d_1(x,-\mu)\\
d_1(x,\mu)-d_2(x,-\mu) &
d_2(x,\mu)-d_1(x,-\mu)
\end{array}
\right]
\\
&+\frac{i}{2\mu}\int_0^x\sigma^2(t)\left[\begin{array}{cc}
e^{i\mu x}-e^{-i\mu x}e^{2i\mu t} &
e^{i\mu x}e^{-2 i\mu t}-e^{-i\mu x} \\
e^{i\mu x}+e^{-i\mu x}e^{2i\mu t} &
e^{i\mu x}e^{-2i\mu t}+
e^{-i\mu x}
\end{array}\right]\,dt
\\
&+\frac{i}{2\mu}
\left[
\begin{array}{cc}
q(x,-\mu) & q(x,\mu) \\
-q(x,-\mu) & q(x,\mu)
\end{array}
\right)
\int_0^x\sigma^2(t)\left(\begin{array}{cc}
1 & e^{-2 i\mu t} \\
-e^{2i\mu t} & -1
\end{array}\right]\,dt,
\end{align*}
\[
 q(x,\mu)=-\int_0^x e^{i\mu(x-2t)} \sigma(t)\,dt,\quad
d_1(x,\mu)=\int_0^x e^{i\mu(x-2t)}\tilde{\sigma}(x,t)\,dt,
\]
\[
d_2(x,\mu)=q(x,\mu)-\int_0^x e^{i(x-2t)} N_\sigma(x,t)\,dt,
\]
where $N_\sigma$ is defined by \eqref{nsigma}.
Moreover, for $\mu\in \Omega_d(\sigma)$,
\begin{equation}\label{phiEst}
\|W_0(x,\mu)\|_{M(C[0,1]}
\le C \rho(\mu),
\end{equation}
where a constant $C>0$  depends only on $d>0$ and $\|\sigma\|_{L_2}$.
\end{theorem}

If in  Theorem
\ref{YSa} we set $\vec{c}=(0,1)^T$, then using the relations
\begin{align}\label{C111}
2\int_0^x \cos(2\mu t)\tilde{\sigma}(x,t)dt
&=\int_0^x\int_0^x \cos(2\mu(t-s))\sigma(s)\sigma(t)\,ds dt
\nonumber\\
&=\left(\int_0^x \cos(2\mu t)\sigma(t)\,dt\right)^2+
\left(\int_0^x \sin(2\mu t)\sigma(t)\,dt\right)^2,\
\end{align}
\begin{align}\label{S111}
\int_0^x \sin(2\mu t)\tilde{\sigma}(x,t)dt
&=\int_0^x\int_0^s\sigma(t)\sigma(s)
\sin(2\mu (s-t)) d t ds\nonumber
\\
&=\int_0^x \sigma(t)\cos(2\mu t)\,dt
\int_0^x\sigma(s)
\sin(2\mu s) ds\nonumber \\ &
-2\int_0^x\int_0^s\sigma(t)\sigma(s)
\cos(2\mu s)\sin(2\mu t) d t ds.
\end{align}
and simple calculations, we get the next corollary.

\begin{corollary}\label{fund}
If $y=y(x,\mu)$ is a solution of \eqref{SE10} with the initial  conditions
\begin{align}\label{conAB}
y(0,\mu)=0,\quad y^{[1]}(0,\mu)=1,
\end{align}
then for $\mu\in \Omega_d(\sigma)$, given by \eqref{domain}, we have
\begin{align}
\mu y(x,\mu)&=\sin(\mu x) +\sin(\mu x)\int_0^x\cos(2\mu t)\sigma(t)dt-\cos(\mu x)\int_0^x\sin(2\mu t)\sigma(t)dt\nonumber\\
&+2\cos(\mu x)\int_0^x\int_0^t\sigma(s)\sigma(t)\cos(2\mu t)\sin(2\mu s) ds dt\nonumber\\
&-\frac{1}{2\mu}\cos(\mu x)\int_0^x\cos(2\mu t)\sigma^2(t)dt
-\frac{1}{2\mu}\sin(\mu x)\int_0^x\sin(2\mu t)\sigma^2(t)dt\nonumber\\
&+\frac{1}{2\mu}\cos(\mu x)\int_0^x \sigma^2(t)dt+  \int_0^x\sin(\mu(x-2t))
N_\sigma(x,t)dt \nonumber\\
&- \cos(\mu x)\int_0^x\sigma(t)\sin(2\mu t)dt\int_0^x\sigma(s)\cos(2\mu s) d s \nonumber\\
&+\frac{1}{2}\sin(\mu x)\Bigg(\left(\int_0^x \cos(2\mu t)\sigma(t)\,dt\right)^2+\left(\int_0^x \sin(2\mu t)\sigma(t)\,dt\right)^2\Bigg)\nonumber\\
&-\frac{1}{2\mu}
\int_0^x\sigma^2(t)dt\int_0^x\cos(2\mu (x-s))\sigma(s)ds
 \nonumber\\
&+\frac{1}{2\mu}\cos(\mu x)\int_0^x\int_0^x\cos(2\mu (t-s))\sigma^2(t)\sigma(s) dt ds
\nonumber\\
&-\frac{1}{2\mu}\sin(\mu x)\int_0^x \int_0^x
\sin (2 \mu(t-s))\sigma^2(t)\sigma(s)dt ds%\nonumber\\
%&
+   \phi_1(x,\mu),\label{y111}
\end{align}
\begin{align}
 y^{[1]}(x,\mu)&=\cos(\mu x) +\cos(\mu x)\int_0^x\cos(2\mu t)\sigma(t)dt+\sin(\mu x)\int_0^x\sin(2\mu t)\sigma(t)dt\nonumber\\
&-2\sin(\mu x)\int_0^x\int_0^t\sigma(s)\sigma(t)\cos(2\mu t)\sin(2\mu s) ds dt\nonumber\\
&+\frac{1}{2\mu}\sin(\mu x)\int_0^x\cos(2\mu t)
\sigma^2(t)dt-\frac{1}{2\mu}\cos(\mu x)\int_0^x\sin(2\mu t)\sigma^2(t)dt\nonumber\\
&-\frac{1}{2\mu}\sin(\mu x)\int_0^x \sigma^2(t)dt +  \int_0^x\cos(\mu(x-2t))
N_\sigma(x,t)dt\nonumber\\
&+ \sin(\mu x)\int_0^x\sigma(t)\sin(2\mu t)dt
\int_0^x\sigma(s)\cos(2\mu s) ds \nonumber\\
&+\frac{1}{2}\cos(\mu x)\Bigg(\left(\int_0^x \cos(2\mu t)\sigma(t)\,dt\right)^2+\left(\int_0^x \sin(2\mu t)\sigma(t)\,dt\right)^2\Bigg)\nonumber\\
&+\frac{1}{2\mu}\int_0^x\sigma^2(t)dt\int_0^x
\sin(\mu(x-2s))\sigma(s)ds
 \nonumber\\
&-\frac{1}{2\mu}\sin(\mu x)\int_0^x\int_0^x
\cos(2\mu(t-s))\sigma^2(t)\sigma(s)dt ds\nonumber\\
&-\frac{1}{2\mu}\cos(\mu x)\int_0^x \int_0^x
\sin(2 \mu (t-s))\sigma^2(t)\sigma(s)dt ds %\nonumber\\
%&
+   \phi_2(x,\mu).\label{y111prim}
\end{align}
where $N_\sigma$ is given by \eqref{nsigma} and
\begin{equation}\label{phiEst}
\|\phi_j(x,\mu)\|_{M(C[0,1])}
\le C \rho(\mu),\quad j=1,2,
\end{equation}
and $\rho$ is defined by \eqref{rhoA}.
\end{corollary}

\section{Spectral problems for Sturm--Liouville equations
with singular potentials
}

In this section we  apply Corollary \ref{fund}
to the spectral problem
\begin{align}\label{specy}
y''(x)&+
q(x)y(x)+\lambda y(x)=0, \ \ \ \ x\in[0,1],\\
y(0)&=0, \ \ y(1)=0.\label{specyw}
\end{align}
where $q=\sigma'$ and $\sigma\in L_2[0,1]$.

Recall that basic facts on asymptotic of eigenvalues and eigenfunctions of the problem \eqref{specy}-\eqref{specyw}
can be found in \cite {S} and \cite{SS}.
 For the reader's convenience and for the matter of comparison with our results,
 we  formulate some of them below.

\begin{theorem}\label{evalue}
If $(\lambda_n)_{n \ge 1}$ are the eigenvalues of the spectral problem \eqref{specy}-\eqref{specyw},
then
\begin{equation}\label{tezaww}
\lambda_n=\mu_n^2,\qquad \mu_n=\pi n + \mu_{0,n}+r_n, \qquad n \in \mathbb N,
\end{equation}
where
\begin{align}\label{tezaww11}
\mu_{0,n}&:=
\int_0^1 \sin(2 \pi n t ) \sigma(t)dt
-2\int_0^1\int_0^t \sigma(t)\sigma(s)\sin (2\pi n s) \cos (2 \pi n t)  ds dt,
\nonumber\\
&-\frac{1}{2\pi n }\int_0^1 (1-\cos(2\pi n t))\sigma^2(t)dt,
\end{align}
and  $(r_n)_{n\ge 1}\in l_1$.
The eigenfunctions $(y_n)_{n \ge 1}$ of the spectral problem \eqref{specy}-\eqref{specyw}  satisfy
\begin{align}
\pi n y_n(x)&= y_{0,n}(x)+\psi_{1,n}(x),\label{tezafw2}
\end{align}
\begin{align}
&y_{0,n}(x)=
\sin(\pi nx)\Bigg(1+\int_0^x \cos(2 \pi ns)  \sigma(s)ds
-\frac{1}{2\pi n}\int_0^x\sin(2\pi nt)\sigma^2(t)dt\Bigg)\nonumber\\
&
+\cos(\pi nx)\Bigg(\mu_{0,n}x -\int_0^x \sin(2 \pi ns ) \sigma(s)ds
+\frac{1}{2\pi n }\int_0^x (1-\cos(2\pi nt))\sigma^2(t)dt
\nonumber\\
&+2\int_0^x\int_0^t \sigma(t)
\sigma(s)\sin (2\pi n s) \cos (2 \pi n t)  ds dt\Bigg),\label{tezafw200}
\end{align}
and
\begin{align}
y_n^{[1]}(x)&=y_{0,n}^{[1]}(x)+\psi_{2,n}(x),\label{tezafw3}
\end{align}
\begin{align}
&y_{0,n}^{[1]}(x)=
\cos(\pi nx)\Bigg(1+\int_0^x \cos(2 \pi ns)  \sigma(s)ds
-\frac{1}{2\pi n}\int_0^x\sin(2\pi nt)\sigma^2(t)dt\Bigg)
\nonumber\\
& -\sin(\pi nx)\Bigg(\mu_{0,n}x-\int_0^x \sin(2 \pi ns ) \sigma(s)ds
+\frac{1}{2\pi n }\int_0^x (1-\cos(2\pi nt)\sigma^2(t)dt
\nonumber\\
&+2\int_0^x\int_0^t \sigma(t)\sigma(s)
\sin (2\pi n s) \cos (2 \pi n t)  ds dt\Bigg),\label{tezafw300}
\end{align}
where
\[
\sup_{x\in[0,1]}\sum_{n=1}^\infty |\psi_{j,n}(x)|<\infty,\quad j=1,2.
\]
\end{theorem}

Our main aim in this section is to provide sharp asymptotic formulas for $(y_n)_{n \ge 1}$.

\begin{theorem}\label{wlasne}
Let $\mu_{0,n}$ be defined by
\eqref{tezaww11} and let  $y_{0,n}(x)$, $y_{0,n}^{[1]}(x), n \ge 1,$ be given by
\eqref{tezafw200}, \eqref{tezafw300}.
The eigenfunctions of the spectral problem \eqref{specy}-\eqref{specyw} admit  the representation
\[
\pi n y_n(x)=y_{0,n}(x)+y_{1,n}(x)+\tilde{\psi}_{1,n}(x),
\]
\[
y_n^{[1]}(x)=y_{0,n}^{[1]}(x)+y_{1,n}^{[1]}(x)+\tilde{\psi}_{2,n}(x),
\]
where\begin{align}
&y_{1,n}(x)=\sin(\pi nx)A_{n}(x)+\cos(\pi nx)B_{n}(x)\label{tezafw}
\end{align}
\begin{align}
&y^{[1]}_{1,n}(x)=
\cos(\pi nx)A_n(x)-\sin(\pi nx)B_n(x)\label{tezafwprim}
\end{align}
and
\begin{align*}
&A_n(x)=
\frac{1}{2}\left(\int_0^x \cos(2\pi n t)\sigma(t)\,dt\right)^2+\frac{1}{2}\left(\int_0^x \sin(2\pi n t)\sigma(t)\,dt\right)^2\nonumber\\
&+ \int_0^x\cos(2\pi nt)
N_\sigma(x,t)dt
\nonumber\\
&-\mu_{0,n}\Bigg(2\int_0^x  \sin(2\pi n t) \sigma(t)tdt
  +
  \int_0^x\sin(2\pi nt)N_\sigma(x,t)tdt\Bigg)\nonumber\\
&+\mu_{0,n}
x \Bigg(\int_0^x \sin(2 \pi nt ) \sigma(t)dt
\nonumber-2\int_0^x\int_0^t \sigma(t)\sigma(s)\sin (2\pi n s) \cos (2 \pi n t)  ds dt\nonumber\\
&+\int_0^x\sin(2\pi n s)\sigma(s)ds\int_0^x\cos(2\pi n t)\sigma(t)
d t+ \int_0^x\sin(2\pi nt)N_\sigma(x,t)dt\Bigg),
\end{align*}
\begin{align*}
&B_n(x)=\mu_{0,n}
x \Bigg(\int_0^x \cos(2\pi n t)\sigma(t)\,dt+\frac{1}{2}\left(\int_0^x \cos(2\pi n t)\sigma(t)\,dt\right)^2\nonumber\\
&+\frac{1}{2}\left(\int_0^x \sin(2\pi n t)\sigma(t)\,dt\right)^2
+ \int_0^x\cos(2\pi nt)
N_\sigma(x,t)dt\Bigg)\nonumber\\
&-\int_0^x\sin(2\pi n s)\sigma(s)ds\int_0^x\cos(2\pi n t)\sigma(t)
d t- \int_0^x\sin(2\pi nt)N_\sigma(x,t)dt\nonumber\\
&+\mu_{0,n}\Bigg(2
\int_0^x  \cos(2\pi n t) \sigma(t)tdt+2\int_0^x\int_0^x\cos\big(2\pi n (s-t)\big)\sigma(s)s\sigma(t)dsdt\nonumber\\
&-  \int_0^x\cos(2\pi nt)N_\sigma(x,t)tdt\Bigg),
\end{align*}
where $N_\sigma$ is given by \eqref{nsigma} and
\[
\sum_{n=1}^\infty\sup_{x\in[0,1]}|\tilde{\psi}_{j,n}(x)|<\infty,\quad
j=1,2.
\]
\end{theorem}

\begin{proof}
Using the relations
\begin{equation*}
\sin(2t\mu_n)-\sin(2\pi n t)
=2\mu_{0,n}t\cos(2\pi nt)
+s_{n}(t),
\end{equation*}
and
\begin{equation*}
\cos(2t\mu_n)-\cos(2\pi n t)
=-2\mu_{0,n}t\sin(2\pi nt)
+c_{n}(t),
\end{equation*}
where
\[
\sum_{n=1}^\infty \sup_{t\in[0,1]}(|s_n(t)|+|c_n(t)|)<\infty,
\]
it suffices to insert the formulas  \eqref{tezaww} and
\eqref{tezaww11} into the identities
\eqref{y111} and \eqref{y111prim},
and to finish the proof by showing that
\begin{equation}\label{ApL1}
\sum_{n=1}^\infty\rho(\mu_n)<\infty,
\end{equation}
where $\rho$ is given by
\eqref{rhoA}.
In what follows we will use a basic formula for eigenvalues
\[\mu_n=\pi n + \tilde{\mu}_n,\] where $(\tilde{\mu}_n)\in l_2$.
To prove \eqref{ApL1},
using  Parseval's identity and
a simple inequality
$|e^{iz}-1|\le |z|e^{d}$, $z\in P_d$,
we infer that
\begin{align}
\sum_{n=1}^\infty
\Big|\int_0^xe^{\pm 2i\mu_nt}\sigma(t)dt\Big|^2
&\le 2\sum_{n=1}^\infty
\left|\int_0^xe^{\pm 2\pi i n t}\sigma(t)dt\right|^2\nonumber\\
&+
2\sum_{n=1}^\infty
\left(\int_0^x|e^{2i\tilde{\mu}_n t}-1||\sigma(t)|dt\right)^2 \nonumber\\
&\le 2
\|\sigma\|^2_{L_2[0,x]}+
4e^{2d}\|\sigma\|_{L_1[0,x]}^2\sum_{n=1}^\infty
|\tilde{\mu}_n|^2\nonumber\\
&\le m<\infty,\label{coeff}
\end{align}
for any $x\in [0,1]$,
where $m:=2(1+
2e^{2d}\|\tilde{\mu}_n\|_{l_2}^2)\|\sigma\|^2_{L_2}$.
It follows from \eqref{coeff} that
\begin{equation}\label{Nart}
\sum_{n=1}^\infty \gamma_0^2(x,\mu)\le  4m,\quad x\in [0,1]
\end{equation}
Finally, by \eqref{Nart},
%if
%\[
%\rho(\mu)=\gamma^2(\mu)+\gamma_1(\mu)+ |\mu|^{-1}
%k_{\sigma^2}(\mu)
%+|\mu|^{-2}
%\]
%then
we have
\begin{align*}
\sum_{n=1}^\infty \rho(\mu_n)&\le
2\int_0^1(1+|\sigma(s)|)\left(\sum_{n=1}^\infty \gamma_0^2(s,\mu_n)\right)\,ds
%+\int_0^1 \sigma_0(s)\left(\sum_{n=1}^\infty\gamma^2_0(s,\mu_n)\right)\,ds
\\
&+2\int_0^1 |\sigma(s)|^2\left(\sum_{n=1}^\infty|\mu_n|^{-1}|\gamma_0(s,\mu_n)\right)\,ds
+\sum_{n=1}^\infty |\mu_n|^{-2}
\\
&\le
2\int_0^1(1+|\sigma(s)|+\|\sigma\|_{L_2}^2)\left(\sum_{n=1}^\infty \gamma_0^2(s,\mu_n)\right)\,ds\\
&+(1+\|\sigma\|_{L_2}^2)\sum_{n=1}^\infty |\mu_n|^{-2}<\infty,
%&\leq 8m+4m\|\sigma_0\|_{L_1}\\
%&+\|\mu_n^{-1}\|_{l_2}\int_0^1 |\sigma(x)|^2\left(\sum_{n=1}^\infty \gamma^2_0(x,\mu_n)\right)^{1/2}\,dx
%+ \|\mu_n^{-2}\|_{l_1}
%\\
%&\le 8m+4m\|\sigma_0\|_{L_1}+2m^{1/2}\|\mu_n^{-2}\|_{l_1} \|\sigma\|^2_{L_2}
%+ \|\mu_n^{-2}\|_{l_1}<\infty.
\end{align*}
thus \eqref{ApL1} follows.
\end{proof}

Note that,  in the same way, we can study a more general Sturm--Liouville equation
\begin{equation}\label{SE0Y}
(p(x)y'(x))'+q(x)y(x)+\mu^2 p(x) y(x)=0,\quad x\in [0,1],\quad
\end{equation}
where $q=u',$ $u\in L_2[0,1],$
and the coefficient $p$ is such that
\[
p\in W_2^1[0,1],\quad
\quad p(x)>0,\quad x\in [0,1].
\]
Indeed,  rewrite  \eqref{SE0Y} as
\begin{equation}\label{SE011}
(y^{[1]}(x))'-\frac{u(x)}{p(x)}y^{[1]}(x)+\frac{u^2(x)}{p(x)}y(x)
+\mu^2 p(x) y(x)=0,\quad x\in [0,1],
\end{equation}
where
\[
y^{[1]}(x):=p(x)y'(x)+u(x)y(x),\qquad y'(x)=\frac{y^{[1]}(x)-u(x)y(x)}{p(x)}.
\]
We say that $y$ is a solution of \eqref{SE011} if
\[
y\in \mathcal{D}:=\{y\in W_2^1[0,1]: y^{[1]}\in W_1^1[0,1]\}
\]
and \eqref{SE011} is satisfied for a.e. $x\in [0,1]$.
We will show that  \eqref{SE011}
can be transformed into the perturbed Dirac system
\eqref{CaychyInt} with appropriate unknown functions $v_1$, $v_2$ and
coefficients $J\in M(L_2[0,1])$ and $P\in M(L_1[0,1])$.

The matrix form of \eqref{SE011} is
 \begin{equation}\label{Znew11X}
\tilde{L}\left(\begin{array} {c}
y \\
y^{[1]}\end{array}
\right):=
\frac{d}{dx}
\left(\begin{array} {c}
y \\
y^{[1]}\end{array}
\right)+M(x,\mu)\left(\begin{array} {c}
y \\
y^{[1]}\end{array}
\right)
=0,
\end{equation}
\begin{equation}
M(x,\mu):=\frac{1}{p}
\left(\begin{array}{cc}
u & -1\\
u^2+\mu^2 p^2 & -u\end{array}\right).
\end{equation}
Define
\[
w(x):=u(x)+p(x)g(x)\in L_2[0,1],\quad
g(x):=\int_0^x \frac{u(t)p'(t)}{p^2(t)}\,dt\in W_2^1[0,1],
\]
\[
\tilde{S}_1(x,\mu)=
\frac{1}{\sqrt{p}}
\left(\begin{array} {cc}
1  & 1 \\
(i\mu-g)p &
-(i\mu+g)p
\end{array}
\right),\]
\[
(\tilde{S}_1(x,\mu))^{-1}
=\frac{1}{2i\mu\sqrt{p}}
\left(\begin{array} {cc}
(i\mu+g)p & 1 \\
(i\mu-g)p &
-1
\end{array}
\right),\quad \mu\not=0.
\]
Since
$2(\sqrt{p})'=p'/\sqrt{p},$ $g'=u p'/p^2,$ we have
\[
\frac{d}{dx}
\tilde{S}_1(x,\mu)
=\frac{p'}{p\sqrt{p}}\left(\begin{array} {cc}
-1 & -1 \\
(i\mu-g)p/2-u&
-(i\mu+g)p/2-u
\end{array}
\right).
\]

For $\mu\not=0$ and $y\in \mathcal{D}$, we introduce a transformation
$y\mapsto \{v_1,v_2\}^T$ by
\[
Y=\tilde{S}_1(x,\mu)V,\quad
Y:=\left(\begin{array} {c}
y \\
y^{[1]}\end{array}
\right),\quad V:=
\left(\begin{array} {c}
v_1 \\
v_2\end{array}
\right).
\]
Then straightforward calculations reveal that
\[
\tilde{L}Y=\tilde{S}_1(\mu)\tilde{L}_0V\qquad
\tilde{L}_0V:=V'+M_0V,
\]
where
\[
M_0=\tilde{S}_1^{-1}(\tilde{S}_1'+M\tilde{S}_1)=
-i\mu J_0-\sigma(x) J-
\frac{i\tau (x)}{2\mu}S,\qquad x\in [0,1].
\]
and
\begin{equation}\label{rfun011}
\sigma(x)=%\frac{u}{p}+g-\frac{p'}{2p}=
\frac{2w-p'}{2p}\in L_2[0,1],
\quad
\tau(x)=\frac{w(p'-w)}{p^2}\in L_1[0,1].
\end{equation}

Summarizing the above, $y\in \mathcal{D}$ is a solution of
\eqref{SE0Y}
if and only if $V=\{v_1,v_2\}^T$ is a solution of
 \eqref{CaychyInt}
with coefficients $\sigma$ and $\tau$ defined by \eqref{rfun011}.
Then one can derive asymptotic formulas  for the fundamental system of solutions
to
\eqref{SE0Y} (as $\mu\to\infty$, $|{\rm Im}\mu|\le c$).
We do not give any details  here since they are quite technical.

\begin{remark}\label{remarc}
We can also  consider the equation
\begin{equation}\label{SE0AY}
(a(x)y'(x))'+q_0(x)y(x)+\mu^2 c(x) y(x)=0,\qquad x\in [0,1],\quad
\end{equation}
where $q_0=u_0',$ $u\in L_2[0,1],$
and the coefficients $a$, $c$ ar such that
\[
a\in W_2^1[0,1],\quad c\in W_2^1[0,1],
\quad a(x)>0,\quad c(x)>0,\quad x\in [0,1].
\]
In this setting,  our results are also applicable, however one needs to employ an additional change of variables.

Let
\[
t=\frac{1}{d}\int_0^x\frac{c^{1/2}(s)}{a^{1/2}(s)}\,ds\in [0,1],\quad
d=\int_0^1\frac{c^{1/2}(s)}{a^{1/2}(s)}\,ds,
z(t)=y(x(t)),\quad t\in [0,1].
\]
Then  from  \eqref{SE0AY} it follows that
\[
\left(p(t)
z'(t)\right)'
+q(t)z(t)+\tilde{\mu}^2 p(t) z(t)=0,\quad t\in [0,1],
\]
where
\[
\tilde{\mu}=d\mu,\quad p(t)=\sqrt{a(x(t))c(x(t))}\in W_2^1[0,1], \quad
q(t)=d^2\frac{\sqrt{a(x(t))}}{\sqrt{c(x(t))}}q_0(x(t)),
\]
and
\[
q(t)=u'(t),\quad u(t)=du_0(x(t))\in L_2[0,1].
\]
This is  exactly the situation described above.
\end{remark}

\section{Appendix}

In this section we include series of results related to the operator $\tilde{T}$, which is used in Section
2, and it is defined by \eqref{Ttilde}. These inequalities are crucial in the proofs concerning the asymptotic
behaviour of the solutions to \eqref{CaychyInt}.

\begin{proposition}\label{4.9}
If $\sigma_j\in L_2[0,1]$, and $F\in M(C(\Delta))$, then
\begin{equation}\label{ind0}
\int_0^xe^{-2i\mu t}(\tilde{T} F)(x,t)dt
=-\int_0^xe^{-2i\mu s}J(s)\int_0^{s}e^{2i\mu \xi}F(s,\xi)d\xi ds.
\end{equation}
Moreover,
\begin{align}\label{wyr2}
\int_0^xe^{-2i\mu t} &(\tilde{T}\tilde{J})(x,t)\,dt \nonumber\\
&=-\int_0^xe^{2i\mu y} \Bigg(\int_y^x J(z)e^{-2i\mu z}\, d z\,
\int_0^{y}J^T(\tau)e^{-2i\mu \tau} d\tau
\Bigg)J^T(y) dy.
\end{align}
\end{proposition}

\begin{proof}
Indeed, we have
\begin{align*}
\int_0^xe^{-2i\mu t}(\tilde{T} F)(x,t)dt
&=-\int_0^xe^{-2i\mu t}\int_t^x J(s)
F(s,s-t)ds dt
\\
&=-\int_0^x J(s) \int_0^s e^{-2i\mu t}F(s,s-t)dt ds
\\&=-\int_0^xe^{-2i\mu s} J(s)\int_0^{s}e^{2i\mu \xi}F(s,\xi)\,d\xi \,ds.
\end{align*}
and \eqref{ind0} is true.
Moreover,
\begin{align*}
\int_0^xe^{-2i\mu t}& (\tilde{T}\tilde{J})(x,t)\,dt
=-\int_0^xe^{-2i\mu t}
 \int_0^{x-t}
 J(t+\xi)
  \int_0^t
  J^T(\tau)J^T(\xi+\tau)
   d\tau d\xi dt
 \\
&=-\int_0^x
 \int_0^{x-\xi}e^{-2i\mu t}
 J(t+\xi)
  \int_0^t
  J^T(\tau)J^T(\xi+\tau)
   d\tau dt d\xi
 \\
%&=-\int_0^x
% \int_\xi^x e^{-2i\mu z}e^{2i\mu \xi}
% J(z)
%  \int_0^{z-\xi}
%  J^T(\tau)J^T(\xi+\tau)
%   d\tau dz d\xi
% \\
&=-\int_0^x
 \int_0^{x-\xi}
  \int_{\xi+\tau}^x
 e^{-2i\mu z}e^{2i\mu \xi}
 J(z)
  J^T(\tau)J^T(\xi+\tau)
 dz  d\tau  d\xi
 \\
%&=-\int_0^x
% \int_0^{x-\tau}
%  \int_{\xi+\tau}^x
% e^{-2i\mu z}e^{2i\mu \xi}
% J(z)
%  J^T(\tau)J^T(\xi+\tau)
% dz d\xi  d\tau
% \\
&=-\int_0^x
 \int_\tau^x
  \int_y^x
 e^{-2i\mu z}e^{2i\mu y}e^{-2i\mu \tau}
 J(z)
 J^T(\tau)J^T(y)
 dz dy  d\tau
% \\
%&=-\int_0^x
% \int_0^y
%  \int_y^x
% e^{-2i\mu z}e^{2i\mu y}e^{-2i\mu \tau}
% J(z)J^T(\tau)J^T(y)
% dz d\tau dy
 \\
&=-\int_0^xe^{2i\mu y} \Bigg(\int_y^x J(z)e^{-2i\mu z}\, d z\,
\int_0^{y}J^T(\tau)e^{-2i\mu \tau} d\tau
\Bigg)J^T(y) dy,
\end{align*}
where  we first changed the order of integration, then passed to new variables $z=t+\xi$ and $y=\tau+\xi$ and changed  the order of integration again.
\end{proof}

Note that if $F\in M(C(\Delta))$ then
\begin{align}\label{Zamec}
\left\|\int_0^x e^{-2tA_\mu} F(x,t)\,dt\right\|_{M(\mathbb{C})}
&\le \left\|\int_0^x e^{-2i\mu t}F(x,t)\,dt\right\|_{M(\mathbb{C})}\nonumber\\&+
\left\|\int_0^x e^{2i\mu t}F(x,t)\,dt\right\|_{M(\mathbb{C})},
\end{align}
for every $x\in [0,1]$.

\begin{lemma}\label{Asimp}
If $\mu\in P_d,$ then
\begin{equation}\label{EstIm0}
\left\|\int_0^x e^{-2tA_\mu}\tilde{J}(x,t)\,dt\right\|_{M(C[0,1])}
\le 2e^{2d}\tilde{l}_0\,\gamma(\mu),
\end{equation}
\begin{equation}\label{EstIm100}
\left\|\int_0^x e^{-2tA_\mu}(\tilde{T}\tilde{J})(x,t)\,dt
\right\|_{M(C[0,1])}
\le 2 e^{4d} \tilde{l}_0^2\gamma(\mu),
\end{equation}
\begin{align}\label{EstIm1}
\left\|\int_0^x e^{-2tA_\mu}(\tilde{T}\tilde{J})(x,t)\,dt\right\|_{M(\mathbb{C})}
&\le 2(l_2+1)e^{2d}
\Big(\gamma(\mu)\gamma_0(x,\mu)\nonumber\\&+\gamma_1(\mu)\Big)
,\quad x\in [0,1],
\end{align}
\begin{equation}\label{EstIm3}
\left\|\int_0^x e^{-2tA_\mu}(\tilde{T}^n \tilde{J})(x,t)\,dt
\right\|_{M(C[0,1])} \le 2e^{2n d}
\frac{l_1^{n-2}}{(n-2)!}
\gamma_2(\mu),\quad n\ge 2.
\end{equation}
\end{lemma}

\begin{proof}
%Changing variables and order of integration,
%analogously as we did in \eqref{ind0}, we obtain
Observe that
\begin{align*}
\left\|\int_0^xe^{-2i\mu t}
\tilde{J}(x,t)dt\right\|_{M(\mathbb{C})}&=
\left\|\int_0^x J(s)e^{-2i\mu s}
\int_0^s J(\xi)e^{2i\mu \xi} d \xi ds\right\|_{M( \mathbb{C})}
\\
&=\left|\int_0^x e^{-2i\mu s}\sigma_1(s)
\int_0^s e^{2i\mu \xi}\sigma_2(\xi) d \xi ds\right|\\&+
\left|\int_0^x e^{-2i\mu s}\sigma_2(s)
\int_0^s e^{2i\mu \xi}\sigma_1(\xi) d \xi ds\right|
\\
&\le  e^{2d}\Big\{\|\sigma_1\|_{L_2}\Big\|\int_0^s e^{2i\mu \xi}\sigma_2(\xi)d\xi \Big\|_{L_2}
\\&+\|\sigma_2\|_{L_2}\Big\|\int_0^s e^{2i\mu \xi}\sigma_1(\xi)d\xi \Big\|_{L_2}\Big\}
\\
&\le  e^{2d}\max\{\|\sigma_1\|_{L_2},\|\sigma_2\|_{L_2}\}\,\gamma(\mu),\quad x\in [0,1].
\end{align*}
Hence, by \eqref{Zamec}, the estimate \eqref{EstIm0} follows.

Next,
in view of \eqref{ind0}, if $\mu\in P_d$, $x\in [0,1]$ and $F\in M(C(\Delta))$, then
\begin{equation}\label{Nes0}
\left\|\int_0^xe^{-2i\mu t}(\tilde{T} F)(x,t)dt\right\|_{M( \mathbb{C})}
\le e^{2d}\int_0^x \left\|J(s)
\int_0^{s}e^{2i\mu \xi}F(s,\xi)d\xi\right\|_{M(\mathbb{C})}\,ds,
\end{equation}
and
\begin{equation}\label{Nes1}
\left\|\int_0^xe^{-2i\mu t}(\tilde{T} F)(x,t)dt\right\|_{M(C[0,1])}
\le e^{2d}l_0
\left\|\int_0^{s}e^{2i\mu \xi}F(s,\xi)d\xi\right\|_{M(C[0,1])}.
\end{equation}
Using \eqref{Nes1} and \eqref{EstIm0}, we obtain that
\begin{align*}
\left\|\int_0^x e^{-2i\mu t}(\tilde{T}\tilde{J})(x,t)\,dt
\right\|_{M(C[0,1])}
&\le e^{2d}l_0
\left\|\int_0^{s}e^{2i\mu \xi}\tilde{J}(s,\xi)\,d\xi\right\|_{M(C[0,1])}\\
&\le e^{4d}\tilde{l}_0^2\gamma(\mu),
\end{align*}
thus, by \eqref{Zamec}, the estimate
\eqref{EstIm100} holds.

By \eqref{Zamec} and the estimate
\[
\left|\int_0^x \sigma_0(s)
\gamma_0(y,\mu)\,dy\right|\le \|\sigma_0\|_{L_2}\left\|\gamma_0(y,\mu)\right\|_{L_2}
\le (l_2+1)\gamma(\mu),
\]
the inequality \eqref{EstIm1} holds if
\begin{equation}\label{EstIm10}
\Bigg\|
\int_0^xe^{-2i\mu t} (\tilde{T}\tilde{J})(x,t)\,dt\Bigg\|_{M(\mathbb{C})}
\le
 e^{2d}\left(
\gamma_1(\mu)
+\gamma_0(x,\mu)
\int_0^x \sigma_0(s)
\gamma_0(y,\mu)\,dy\right).
\end{equation}
On the other hand, using \eqref{wyr2}, we have
\begin{align*}
\Bigg\|
\int_0^x&e^{-2i\mu t} (\tilde{T}\tilde{J})(x,t)\,dt\Bigg\|_{M(\mathbb{C})}\\
&\le
e^{2d}
\int_0^x\left\|\int_{y}^x e^{-2i\mu z}J(z) d z \int_0^y e^{-2i\mu \tau}J^T(\tau)\,d\tau J^T(y)\right\|_{M_2({\rm C})}
dy
\\
&\le e^{2d}
\int_0^x \sigma_0(y)
\left\|\int_{y}^x e^{-2i\mu z}J(z) d z\right\|_{M(\mathbb{C})}
\left\| \int_0^y e^{-2i\mu \tau}J(\tau)\,d\tau \right\|_{M(\mathbb{C})}
dy
\\
&\le e^{2d}
\int_0^x \sigma_0(y)
\left\| \int_0^y e^{-2i\mu \tau}J(\tau)\,d\tau \right\|^2_{M(\mathbb{C})}
dy
\\
&+e^{2d}\left\|\int_0^x e^{-2i\mu z}J(z) d z\right\|_{M(\mathbb{C})}
\int_0^x \sigma_0(s)
\left\| \int_0^y e^{-2i\mu \tau}J(\tau)\,d\tau \right\|_{M(\mathbb{C})}
dy
\\
&\le e^{2d}
\gamma_1(\mu)
+e^{2d}\gamma_0(x,\mu)
\int_0^x \sigma_0(s)
\gamma_0(y,\mu)\,dy,
\end{align*}
and \eqref{EstIm10} follows.

To prove \eqref{EstIm3} it suffices to show that
for all $n\ge 2$ and any $x\in [0,1],$
\begin{align}\label{EstIm30}
\left\|\int_0^x e^{-2i\mu t}(\tilde{T}^n \tilde{J})(x,t)\,dt\right\|_{M(\mathbb{C})}
\le \frac{e^{2n d}}{(n-2)!}
\left(\int_0^x \sigma_0(s)\,ds\right)^{n-2}
\gamma_2(\mu).
\end{align}
We prove \eqref{EstIm30}  by induction.
Using \eqref{Nes0}
for $F=\tilde{T}\tilde{J}$
and  \eqref{EstIm10},
%\eqref{ga2},
we note that
\begin{align*}
 \Big\|\int_0^x e^{-2i\mu t}(\tilde{T}^2\tilde{J})&(x,t)\,dt\Big\|_{M(\mathbb{C})}
\le e^{2d}\int_0^x \sigma_0(s)
\left\|\int_0^{s}e^{2i\mu \xi}(\tilde{T}\tilde{J})(s,\xi)d\xi\right\|_{M(\mathbb{C})}\, ds
\\
&\le e^{4d}\int_0^x |\sigma_0(s)
\left(\gamma_0(s,\mu)\int_0^s\sigma_0(y)
\gamma_0(y,\mu)\,dy+\gamma_1(\mu)\right)\, ds
\\
&\le e^{4d}\int_0^x \sigma_0(s)\gamma_0(s,\mu)
\int_0^{s}\sigma_0(y)
\gamma_0(y,\mu)\,dy\, ds
+e^{4d}l_1\gamma_1(\mu)
\\
&\le
 e^{4d}\frac{\left(\int_0^x \sigma_0(s)\gamma_0(s,\mu)\,ds\right)^2}{2}
+e^{4d}l_1\gamma_1(\mu)
\\
&\le e^{4d}l_2\gamma^2(\mu)+e^{4d}l_1\gamma_1(\mu).
\end{align*}
Therefore, \eqref{EstIm30} holds for $n=2$.

Suppose now that \eqref{EstIm30}
is true for some $n\ge 2$.
Then, once again using \eqref{Nes0}, we have
\begin{align*}
\Big\|\int_0^x e^{-2i\mu t}(\tilde{T}^{n+1}\tilde{J})&(x,t)\,dt\Big\|_{M(\mathbb{C})}\\
&\le e^{2d} \int_0^x \sigma_0(s)
\left\|\int_0^{s}e^{2i\mu \xi}(\tilde{T}^n\tilde{J})(s,\xi)d\xi
\right\|_{M(\mathbb{C})} ds
\\
&\le \frac{e^{2(n+1) d}}{(n-2)!}\gamma_2(\mu)\int_0^x\sigma_0(s)
 \left(\int_0^s\sigma_0(\tau)\,d\tau\right)^{n-2}\,ds
\\
&=\frac{e^{2(n+1)d}\gamma_2(\mu)}{(n-1)!}
 \left(\int_0^x \sigma_0(\tau)\,d\tau\right)^{n-1},\quad x\in [0,1],
\end{align*}
so \eqref{EstIm30}
holds also for  $n+1,$ and  the proof
of \eqref{EstIm30} is finished.
\end{proof}

\bibliographystyle{plain}
\bibliography{bibi2}

\end{document}